\documentclass{article}[12pt]
\usepackage{amsfonts}
\usepackage{amsmath}
\usepackage{graphicx}
\usepackage{amssymb}
\usepackage{epstopdf}
\usepackage{hyperref}
\usepackage[all]{xy}
\usepackage[english]{babel}

\newtheorem{thm}{Theorem}    
\newtheorem{lem}{Lemma}[section]      
\newtheorem{prop}[lem]{Proposition}
\newtheorem{defn}[lem]{Definition}    
\newtheorem{remark}[lem]{Remark}

\title{On computation of morphism spaces and a direct limit of the bordered Floer homology of knot complements}

\author{Jaepil Lee}

\begin{document}

\maketitle

\begin{abstract}
In the bordered Floer theory, gluing thickened torus of positive meridional Dehn twist to the boundary of a knot complement result in the knot complement of increased framing. For a fixed knot $K$, we construct a direct system of positively framed knot complements and study the direct limit. We also study the morphism space between two type-$DD$ modules, and derive type-$DA$ morphisms from $DD$ morphisms to derive the direct system maps. In addition, we introduce a direct limit invariant from the direct system which can detect non-unstable chains in the type-$D$ module of a knot complement, if the type-$D$ modules of the direct system are obtained by algorithm of Lipshitz, Ozsv\'ath and Thurston.
\end{abstract}

\section{Introduction}
Since the Heegaard Floer homology made a breakthrough in low dimensional topology society, there have been many interesting variants of Heegaard Floer homology recently invented. Especially, the bordered Heegaard Floer homology was the one of the most interesting variant for three manifolds with boundaries, which is equipped with the $\mathcal{A}_{\infty}$ module structure or type-$D$ structure. The Pairing Theorem introduced in \cite{LOT08} glues these two modules and the resulting complex is homotopic to the Heegaard Floer homology of closed three manifold. \\

In \cite{LOT08}, the complete description to recover the type-$D$ module of bordered Heegaard Floer homology of knot complement with framing numbers from classical knot Floer homology was given. Roughly speaking, the type-$D$  module consists of two parts; a stable chain and a unstable chain. The stable chain is the part which can be directly recovered from knot Floer homology and independent from the framing number. On the other hand, the homotopy type of the unstable chain is dependent on the framing number. \\

In this paper, we study to distinguish the stable and unstable chain by constructing a direct system $(A_i, \phi_{ij})$ where $A_i$ is a type-$D$ module of knot complement with sufficiently large positive framing $i_0 +i$, and the map $\phi_{ij} : A_i \rightarrow A_{i+j} $ is a map between type-$D$ modules. 
We will define an alternate, coordinate free definition of the unstable complex by using the direct limit of the direct system $(A_i, \phi_{ij})$. In order to do so, we will first need to find all homotopically nontrivial maps between the identity module $\widehat{CFDA}(\mathbb{I})$ to the framing-increasing module $\widehat{CFDA}(\tau_{m+1})$, and glue it to the type-$D$ module of knot complement $\widehat{CFD}(S^3 \backslash \textrm{nbd}(K))$. By doing so, we obtain a map between two type-$D$ modules of the same knot complement whose framing difference is one. However, instead of studying the morphism space $\textrm{Mor}^{\mathcal{A}(-T^2)} (\widehat{CFDA}(\mathbb{I}), \widehat{CFDA}({\tau_{m+1}}) )_{\mathcal{A}(T^2)} $, which has fairly sophisticated structure, we study the morphism spaces between two type-$DD$ bimodules with less structures underneath. Those two morphism spaces have the same homotopy type because of the equivalence of categories proved in \cite{LOT10}. \\

This paper is structured in the following way. In section 2, we quickly review the materials used in this paper. In particular, by using a similar analysis used in \cite{LOT11}, we identify the morphism space between two type-$DD$ modules to a Heegaard Floer homology of some closed Heegaard diagram obtained by gluing two Heegaard diagrams with two boundary components, with a ``full boundary twist'' inserted within (this will be explained in section 2). This proves the following theorem.
\begin{thm}
Let $Y_1 = Y( \mathcal{H}_1)$ and $Y_2 = Y(\mathcal{H}_2)$ be three manifolds with two boundary components. Then the morphism space $\mathrm{Mor}^{\mathcal{A}(-T^2), \mathcal{A}(-T^2)} (\widehat{CFDD}(Y_1), \widehat{CFDD}(Y_2))$ is isomorphic (up to homotopy) to
\begin{displaymath}
\widehat{CFAA} (\tau^{-1}_{\partial} ( Y(\mathcal{H}_2) ) \boxtimes_{\mathcal{A}(T^2), \mathcal{A}(-T^2)} \widehat{CFDD} ( -Y(\mathcal{H}_1) ),
\end{displaymath}
where $\tau^{-1}_{\partial} (Y (\mathcal{H}_2))$ implies a Heegaard diagram obtained from $\mathcal{H}_2$ with the full negative boundary twist diagram attached in the one end.
\end{thm}
Note that the same consequence can be derived from the result in \cite{LOT11}, Corollary 11.\\

Sections 3 and 4 study the explicit computation of the morphism space and compute the homology $H_* (\textrm{Mor}^{\mathcal{A}(-T^2), \mathcal{A}(-T^2)} (\widehat{CFDD}(\mathbb{I}), \widehat{CFDD}({\tau_{m+1}}) ) )$, and convert type-$DD$ morphisms to $DA$ morphisms by tensoring with the identity type-$AA$ module. In section 5, the direct system and direct limit are introduced and computes the direct limits of three typical chains appearing in $\widehat{CFD}(S^3 \backslash \textrm{nbd}(K) )$. In particular, let $f'$ be any morphism in $\mathrm{Mor}^{\mathcal{A}(-T^2)} (\widehat{CFDA}(\mathbb{I}), \widehat{CFDA}(\tau_{m+1}))_{\mathcal{A}(T^2)}$, and ${}^{\mathcal{A}(-T^2)} A$ be a type-$D$ module. Then $f=f' \boxtimes \mathrm{id}_A : A \rightarrow \widehat{CFDD}(\tau_{m+1}) \boxtimes A$. Moreover, by composing $f$ multiple times we get a map 
\begin{displaymath}
\phi_{ij} := \underbrace{f \circ \cdots \circ f}_{j-i \textrm{ times}} : A_i \rightarrow A_j, 
\end{displaymath}
where $A_i := \underbrace{\widehat{CFDA}( \tau_{m+1} ) \boxtimes \cdots \boxtimes \widehat{CFDA}( \tau_{m+1} ) }_{i \textrm{ times}} \boxtimes A \cong \widehat{CFDA}( \tau_{m+i} ) \boxtimes A$. \\
We consider the direct system $(A_i, \phi_{ij})$ constructed in this manner and compute the direct limit $\varinjlim A_i$. Especially we compute the direct system of three typical types of chains that appear in the type-$D$ modules of knot complements(with positive framing), by changing $A$. \\

Next we consider a direct system $(A_i, \phi_{ij})_{i=0}$, where $A_0$ is a type-$D$ module of a knot $K$ complement with sufficiently large \emph{negative} framing. Such a type-$D$ module can be derived from a \emph{reduced} model(whose differential strictly drops the Alexander of $U$-filtration, see Definition 11.25 of~\cite{LOT08}) of a chain complex $CFK^-(K)$. Each nontrivial differential of the reduced $CFK^-(K)$ that drops either the Alexander filtration or the $U$-filtration is converted to a certain sequence of chains in the type-$D$ module $\widehat{CFD}(S^3 \backslash K)$ of the knot $K$ complement. In Theorem 11.27 of~\cite{LOT08}, they introduce a \emph{unstable chain} in $\widehat{CFD}(S^3 \backslash K)$, whose structure is dependent on the framing data of the knot complement. Roughly, in $\widehat{CFD}(S^3 \backslash K)$, there are two different types of chains - the unstable chain and non-unstable chain. The structure of those non-unstable chains is not affected by the framing and solely depending on the structure of $CFK^-(K)$. However, these chains are only defined only by using a special basis: under the assumption that $\iota_0 \widehat{CFD}(S^3 \backslash K) $ is reduced, the chain is only defined between horizontally of vertically simplified basis of $\iota_0 \widehat{CFD}(S^3 \backslash K)$. See Definition 11.23 of \cite{LOT08}. \\ 

In section 6, we introduce a direct limit invariant that detects non-unstable chains only.

\begin{thm}
Let $A_i$ be a type-$D$ module of a knot $K$ complement with framing $i_0 +i$, where $i$ is any nonnegative integer and $i_0$ is a sufficiently large negative integer. Also assume that the direct limit of a direct system $(A_i, \phi_{ij})_{i=0}^{\infty}$ is nontrivial. For a natural inclusion map $\iota : A_0 \rightarrow \varinjlim A_i$, the homotopy type of $\iota (A_0) \subset \varinjlim A_i$ is a knot invariant. 
Moreover, assume that each $A_i$ is the type-$D$ module of a knot $K$ complement with framing $i_0 +i$, such that $A_i$ is exactly the same complex, given by the algorithm in the Theorem 11.26 of~\cite{LOT08}. Then $\iota (A_0)$ is homotopic to the complex $D$ such that
\begin{itemize}
  \item $D = \iota_0 D \bigoplus \iota_1 D$ as a vector space decomposition and $\iota_0 D = CFK^-(K)$ is  reduced.
  \item Let $\{ \xi_n \}$ be a vertically simplified basis of $\iota_0 D = CFK^-(K)$. Then there is a sequence $\{ \kappa_k^n \} \subset \iota_1 A_i$ where
  \begin{displaymath}
    \xymatrix{
\xi_n \ar[r]^{\rho_1} & \kappa^n_1 & \cdots \ar[l]_{\rho_{23}} & \kappa^n_k \ar[l]_{\rho_{23}} & \kappa^n_{k+1} \ar[l]_{\rho_{23}} & \cdots \ar[l]_{\rho_{23}} & \kappa^n_l \ar[l]_{\rho_{23}} & \xi_{n+1} \ar[l]_{\rho_{123}}.
}
  \end{displaymath}
  \item let $\{ \eta_n \}$ be a horizontally simplified basis of $\iota_0 D = CFK^-(K)$. Then there is a sequence $\{ \lambda_k^n \} \subset \iota_1 A_i$ where
  \begin{displaymath}
    \xymatrix{
\eta_n \ar[r]^{\rho_3} & \lambda^n_1 \ar[r]^{\rho_{23}} & \cdots \ar[r]^{\rho_{23}} & \lambda^n_k \ar[r]^{\rho_{23}} & \lambda^n_{k+1} \ar[r]^{\rho_{23}} & \cdots \ar[r]^{\rho_{23}} & \lambda^n_l \ar[r]^{\rho_2} & \eta_{n+1}.
}
  \end{displaymath}
  \item As a $\mathbb{F}_2$ vector space, $\iota_1 D$ has no other element.  
\end{itemize}
\end{thm}

Note that the module $D$ does not have the unstable chain. This theorem also gives a coordinate free definition of the ``non-unstable chain'' of the type-$D$ module of a knot complement.  

\subsection{Further Questions}
This paper was originally motivated in the program of finding type-$DD$ module of 2-link components from the Ozsv\'ath-Sz\'ab\`o link invariant $CFL^{\infty}$ ~\cite{OS05}. The idea of taking direct limit of knot complement by gluing $T^2 \times I$ was inspired by work of John B. Etnyre, David Shea Vela-Vick and Ruman Zarev \cite{EVZ14}. There was an attempt to find type-$DD$ module of $(2,2n)$ torus link complement \cite{Lee13}, but the computation was unable to distinguish the generators dependent on $CFL^{\infty}$ and generators dependent on the framings of each link component. These are some questions that naturally arise. \\

\textbf{Question 1.} \emph{Let $A_i(K)$ be the type-$D$ module of knot $K$ complement. The maps given in section 4 is the only natural maps from $A_i$ to $A_{i+1}$; that is, these maps preserves the skein relation.} \\

According to \cite{LOT11}, there are more morphisms in $\mathrm{Mor}(A_i, A_{i+1}) $ than $\mathrm{Mor} ( \widehat{CFDA}(\mathbb{I}), \widehat{CFDA}(\tau_{m+1}) ) $, since the morphism space $\textrm{Mor}(A_i,A_{i+1})$ has the same homotopy type as the closed Heegaard Floer homology of the diagram of glued two knot complements. However, not every morphisms are natural in the sense of respecting the skein relation. \\

\textbf{Question 2.} \emph{Are we able to identify the ``stable chain'' of the type-$DD$ module of a 2-link complement by studying the $G'(\mathcal{Z})$-$G'(\mathcal{Z})$ grading of the bimodule?} \\
  
For a 2-link in $S^3$, we may comupute the type-$DD$ left-left bimodule over the strands algebras $\mathcal{A}(- \mathbb{T}^2)$. For the type-$DD$ modules of $(2,2n)$-torus links have been computed in~\cite{Lee13}. By applying same trick to the type-$DD$ module, we hope to recover the ``invariant'' part of the module regardless of the framing.

\subsection{Acknowledgement}
The author thanks to Robert Lipshitz for the helpful discussions about the morphism spaces in the category of the type-$D$ module and assorted bimodules. Especially referring to Tova Brown's thesis~\cite{Br11} for the geometric interpretation of the morphism spaces.
 
\section{Background}
\subsection{Algebraic Preliminaries}
In this subsection we will quickly recall the algebraic tools that will be used in the rest of this paper.\\

Let $\mathcal{A}$ be a $dg$ algbera with differential $\mu_1$ and associative multiplication $\mu_2$. A \emph{type-$D$} module, or \emph{type-$D$} structure ${}^{\mathcal{A}} M$ is a $\mathbf{k}$-module $M$ with left $\mathcal{A}$ action, equipped with a map $\delta^1 : M \rightarrow \mathcal{A} \otimes M$ satisfying
\begin{displaymath}
(\mu_2 \otimes \mathbb{I}_M) \circ ( \mathbb{I}_{\mathcal{A}} \otimes \delta^1 ) \circ \delta^1 + (\mu_1 \otimes M) \circ \delta^1 : M \rightarrow \mathcal{A} \otimes M
\end{displaymath} 
vanishes.\\

The dual of type-$D$ module is constructed as follows. Let $\overline{M} = \mathrm{Hom}_{\mathbf{k}} (M, \mathbf{k})$. The structure map $\delta^1$ of a type-$D$ structure can be interpreted as an element in $\overline{M} \otimes \mathcal{A} \otimes M$, which can also be interpreted as
\begin{displaymath}
\overline{\delta}^1 : \overline{M} \rightarrow \overline{M} \otimes \mathcal{A}.
\end{displaymath}
Then it is an easy exercise to show that $(\overline{M}, \overline{\delta}^1 )$ is a right type-$D$ module. \\

For type-$D$ modules ${}^{\mathcal{A}} M$ and ${}^{\mathcal{A}} N$ with at least one of $M$ or $N$ begin a  finite dimensional, the morphism space $\mathrm{Mor}^{\mathcal{A}}(M,N)$ is isomorphic to $\overline{M} \boxtimes \mathcal{A} \boxtimes N$ as a chain complex (Proposition 2.7,~\cite{LOT11}).\\  

We recall the definition of a \emph{Hochschild complex} (Definition 2.3.41, \cite{LOT10}). A Hochschild complex $CH( {}_{\mathcal{A}} M_{\mathcal{A}} )$ of $\mathcal{A}_{\infty}$-algebra $\mathcal{A}$ and $\mathcal{A}_{\infty}$-bimodule ${}_\mathcal{A} M_{\mathcal{A}}$ is defined as follows. $CH_n ( M )$ is $\mathbf{k}$-vector space, quotient of $M \otimes_{\mathbf{k}} \underbrace{ \mathcal{A} \otimes_{\mathbf{k}} \cdots \otimes_{\mathbf{k}} \mathcal{A} }_{n \textrm{-times} }$ by the relation
\begin{displaymath}
  e \cdot x \otimes a_1 \otimes \cdots \otimes a_n = x \otimes a_1 \otimes \cdots \otimes a_n \cdot e
\end{displaymath}
where $e \in \mathbf{k}$. As a vector space, 
\begin{displaymath}
  CH( {}_{\mathcal{A}} M_{\mathcal{A}} ) = \bigoplus_{n=0}^{\infty} CH_n ( M ).
\end{displaymath}
The differential $D$ on $ CH( {}_{\mathcal{A}} M_{\mathcal{A}} )$ is defined to be
\begin{eqnarray*}
D(x \otimes a_1 \otimes \cdots \otimes a_l) 
& = & \sum_{m+n \leq l} m_{m,1,n}(a_{l-m+1}, \cdots, a_l, x, a_1, \cdots, a_n) \otimes a_{n+1} \otimes \cdots \otimes a_{l-m} \\
& & + \sum_{1 \leq m < n \leq l} x \otimes a_1 \otimes \cdots \otimes \mu_{n-m} (a_m, \cdots, a_{n-1} ) \otimes a_n \otimes \cdots \otimes a_l.
\end{eqnarray*}
Sometimes we drop subscripts $\mathcal{A}$ from the notation if the algebra is clear from its context. \\

There is a more convenient way to understand the Hochschild complex. The summand of Hochschild complex $CH_n (M)$ is generated by $n$ elements from $\mathcal{A}$ and one element from $M$ arranged on a circle. Then, the differential $D$ can be understood as choosing any $i \leq n$ consecutive elements from the circle and apply $i$th order multiplication on the chosen elements. $D^2=0$ can be easily proved from $\mathcal{A}_{\infty}$ relations. \\

We will denote $HH ( M )$ the homology of $CH( M )$.

\subsection{Review on bordered Floer theory}
We will recall the important results from \cite{LOT08}, \cite{LOT10} and \cite{LOT11}, that will be used in the rest of this paper. \\

A \emph{closed Heegaard diagram} is $(\Sigma_g, \boldsymbol{\alpha}, \boldsymbol{\beta} )$, where $\Sigma_g$ is a closed surface with genus $g$ and $\boldsymbol{\alpha} = \{ \alpha_1, \cdots \alpha_g \}$ and $\boldsymbol{\beta} = \{ \beta_1, \cdots, \beta_g \}$. Each of these $g$ sets of circles are non-intersecting circles on $\Sigma_g$, specifying attaching circles of $\Sigma_g \times [0,1]$. A \emph{bordered Heegaard diagram} is $(\overline{\Sigma}_g, \overline{ \boldsymbol{\alpha} }, \boldsymbol{\beta} )$, where $\overline{\Sigma}_g$ is a genus-$g$ surface with a single puncture. $\boldsymbol{\beta}$ is a set of $g$ attaching circles, but $\overline{ \boldsymbol{ \alpha } }$ consists of $g-k$ circles $\boldsymbol{\alpha}^c$ and $2k$ arcs $\boldsymbol{\alpha}^a$ whose boundaries are on $\partial \overline{\Sigma}_g$. \\

A \emph{matched circle} is $\mathcal{Z} := (Z, \mathbf{a}, M)$. $Z$ is an oriented circle, $\mathbf{a}$ is a set of $4k$ points on $Z$, and $M : \mathbf{a} \rightarrow \mathbf{a}$ is a fixed-point free involution. The matching determines a parametrization of surface with genus $k$, and the surface is denoted $F (\mathcal{Z} )$. There is a distinguished disk $D$ on the surface whose boundary is identified to $Z$, sometimes called a \emph{preferred disk}. A \emph{pointed matched circle} is a matched circle with a point $z$ on $Z$ away from $\mathbf{a}$. In this paper, we will be only interested in the torus boundary case (i.e, $k=1$). \\

A \emph{bordered three manifold} is $(Y, \Delta, z, \psi)$, where $Y$ is a three manifold with boundary, $\Delta$ is a disk in $\partial Y$, $z$ is a point on $\partial \Delta$, and $\psi$ is a homeomorphism such that
\begin{displaymath}
\psi : ( F( \mathcal{Z} ), D, z) \rightarrow ( \partial Y, \Delta, z),
\end{displaymath}
so that it gives a parametrization data of the boundary surface. \\

A bordered Floer homology package associates pointed matched circle to a $dg$ algebra $\mathcal{A} ( \mathcal{Z} )$. In general $\mathcal{A} ( \mathcal{Z} )$ has nontrivial differential, but in the torus boundary case the differential is trivial. The torus algebra is a strands algebra with $k=1$, and it has a relatively simply structure; a reader can find its explicit description in the Chapter 11 of \cite{LOT08}. \\

From a three manifold $Y$ with boundary parametrized by $F (\mathcal{Z})$, a left $dg$ $\mathcal{A} ( - \mathcal{Z} )$ module $\widehat{CFD}(Y)$ and a right $\mathcal{A}_{\infty}$ module $\widehat{CFA}(Y)$ can be defined. The negative sign of the $ - \mathcal{Z}$ means the boundary has the opposite orientation from the induced orientation. These two modules are defined from the Heegaard diagram of the three manifold $Y$, and they are well defined up to homotopy. \\

The \emph{pairing theorem} enables the computation of the classical Heegaard Floer homology of the closed three manifold. If $\partial Y_1 = F (\mathcal{Z}) = -\partial Y_2$, then 
\begin{displaymath}
\widehat{CF}(Y_1 \cup_{\partial} Y_2) \cong \widehat{CFA}(Y_1) \tilde{\otimes}_{\mathcal{A} (\mathcal{Z}) } \widehat{CFD}(Y_2),
\end{displaymath} 
where $\tilde{\otimes}$ means the derived tensor product. \\

We can generalize the theory to three manifolds with two boundary components. 
\begin{defn}
  A \emph{doubly bordered Heegaard diagram} ${}^{\alpha} \mathcal{H}^{\alpha}$ with two boundary
  components is a quadruple $(\overline{\Sigma}_g, \boldsymbol{\alpha}, \boldsymbol{\beta},
  \mathbf{z})$ satisfying:
  \begin{itemize}
    \item $\overline{\Sigma}_g$ is a compact, genus $g$ surface with 2 boundary
    components $\partial_L \overline{\Sigma}_g$ and $\partial_R
    \overline{\Sigma}_g$.
    \item $\boldsymbol{\beta}$ is a $g$-tuple of pairwise disjoint
    curves in the interior of $\overline{\Sigma}_g$.
    \item $\boldsymbol{\alpha} = \{ \boldsymbol{\alpha}^{a,L}=\{
    \alpha_1^{a,L}, \cdots, \alpha_{2l}^{a,L} \}, \
    \boldsymbol{\alpha}^{a,R}=\{\alpha_1^{a,R}, \cdots, \alpha_{2r}^{a,R} \}, \ \boldsymbol{\alpha}^{c}=\{\alpha_1^c, \cdots, \alpha_{g-l-r}^c \}
    \}$, is a collection of pairwise disjoint embedded arcs with
    boundary on $\partial_L \overline{\Sigma}_g$ (the
    $\alpha_i^{a,L}$), arcs with boundary on $\partial_R \overline{\Sigma}_g$ (the
    $\alpha_i^{a,R}$), and circles (the $\alpha_i^c$) in the
    interior of $\overline{\Sigma}_g$.
    \item $\mathbf{z}$ is a path in $\overline{\Sigma}_g \backslash
    (\boldsymbol{\alpha} \cup \boldsymbol{\beta})$ between $\partial_L \overline{\Sigma}_g$ and $\partial_R \overline{\Sigma}_g$.
  \end{itemize}
\label{def:doubly}
\end{defn}
Attaching three-dimensional two handles on $ \overline{\Sigma}_g \times [0,1]$ along $\alpha$- and $\beta$- circles will result the three manifolds with two boundary components. $\alpha$-arcs will give the parametrization of the respective boundaries.\\

We associate various bimodules three manifold $Y$ with two boundary component $F(\mathcal{Z}_1$ and $F(\mathcal{Z}_2)$.  
\begin{itemize}
\item $\widehat{CFDD}(Y)$ is a $dg$ bimodule with two left $\mathcal{A} (-\mathcal{Z}_1), \mathcal{A} (-\mathcal{Z}_2)$ actions.
\item $\widehat{CFDA}(Y)$ is an $\mathcal{A}_{\infty}$ bimodule with a left $\mathcal{A} (-\mathcal{Z}_1)$ action and a right $\mathcal{A} (\mathcal{Z}_2)$ action.
\item $\widehat{CFAA}(Y)$ is an $\mathcal{A}_{\infty}$ bimodule with two right $\mathcal{A} (\mathcal{Z}_1), \mathcal{A} (\mathcal{Z}_2)$ actions. 
\end{itemize}
Each of these modules are well defined up to homotopy, and these modules satisfy pairing theorems as well.\\

\subsection{$\alpha$- $\beta$- Bordered Heegaard Diagram}
A three manifold with boundary can be achieved in slightly different ways. In~\cite{LOT08}, the Heegaard diagram only used $2k$ $\alpha$-arcs to parametrize its boundary, not $\beta$-arcs. Since in this paper we will be using both $\alpha$- and $\beta$- arcs for parametrization, we will give additional decoration to the pointed matched circle. First, here we repeat the Definition 3.1 from~\cite{LOT11}.
\begin{defn}
A \emph{decorated pointed matched circle} consists of the following data.
  \begin{itemize}
    \item a circle $Z$.
    \item a decomposition of $Z$ into two closed oriented intervals, $Z_{\alpha}$ and $Z_{\beta}$, whose intersection consists of two points, and $Z_{\alpha}$ and $Z_{\beta}$ oriented opposingly;
    \item a collection of $4k$ points $\mathbf{p} =\{ p_1, \cdots, p_{4k} \}$ in $Z$, so that either $\mathbf{p} \subset Z_{\alpha}$ or $\mathbf{p} \subset Z_{\beta}$.
    \item a fixed-point-free involution $M$ on the points in $\mathbf{p}$; and
    \item a decoration by the letter $\alpha$ or the letter $\beta$, which indicates whether the points $\mathbf{p}$ lie in $Z_{\alpha}$ or $Z_{\beta}$.
  \end{itemize}
\end{defn}
It is easy to see $(Z = Z_{\alpha} \cup Z_{\beta}, \mathbf{p}, M)$ is a matched circle. We will call the matched circle $\mathcal{Z}^{\alpha}$ or $\mathcal{Z}^{\beta}$ depending on the decoration. Then construction of $F ( \mathcal{Z}^{\epsilon} )$ ($\epsilon$ is $\alpha$ or $\beta$) is parallel, but for the further use a slightly different construction will be introduced. \\

Consider a oriented disk $D_0$ with boundary $Z$, so that the orientation of $\partial D_0$ agrees to $Z_{\alpha}$. Let $\{ h, t\} = Z_{\alpha} \cap Z_{\beta}$, and attach one handle $s$ along $h$ and $t$. Finally, attach two disks $D_{\alpha}$ and $D_{\beta}$ to the remaining boundary respecting the decoration. The resulting surface $F (\mathcal{Z}^{\epsilon})$ has a preferred embedded disk, $D_{\alpha} \cup s \cup D_{\beta}$. \\  

A Heegaard diagram $\mathcal{H}$ representing $Y ( \mathcal{H} )$ has the boundary $ \partial Y$ parametrized by $\alpha$-arcs. Then the boundary of $\mathcal{H}$ is a matched circle $\mathcal{Z}^{\alpha}$. From now on, we will call a Heegaard diagram an \emph{$\alpha$-bordered Heegaard diagram} which will be denoted $\mathcal{H}^{\alpha}$. Of course, there is a symmetric construction namely \emph{$\beta$- bordered Heegaard diagram} which uses $\beta$ curves to encode its parametrization information. As the name suggests, its construction is parallel to an $\alpha$-bordered Heegaard diagram. The only difference is the $\beta$- bordered Heegaard diagram uses $2k$ $\beta$ arcs to parametrize its boundary.  For a given $\alpha$-bordered Heegaard diagram $\mathcal{H}^{\alpha}$, one can associate a $\beta$-bordered Heegaard diagram $\overline{\mathcal{H}}^{\beta}$ obtained by interchanging the labels of $\alpha$ and $\beta$ curves. In \cite{LOT11}, it is proved that
\begin{displaymath}
\widehat{CFD} (\overline{\mathcal{H}}^{\beta} )^{\mathcal{A} (- \mathcal{Z})}  \cong {}^{\mathcal{A} (- \mathcal{Z})}  \overline{\widehat{CFD} (\mathcal{H}^{\alpha}) }.
\end{displaymath}
A similar set up exists for a doubly bordered Heegaard diagram. The doubly bordered Heegaard diagram introduced in~\ref{def:doubly} used two sets of $\alpha$-arcs for parametrization. From now on, this doubly bordered Heegaard diagram will be called an \emph{$\alpha$-$\alpha$ bordered Heegaard diagram}, denoted ${}^{\alpha} \mathcal{H}^{\alpha}$. If two boundaries are parametrized by two sets of $\beta$-arcs, it will be called a \emph{$\beta$-$\beta$ bordered Heegaard diagram} and denoted ${}^{\beta} \mathcal{H}^{\beta}$. Likewise there exists an \emph{$\alpha$-$\beta$ bordered diagram} ${}^{\alpha} \mathcal{H}^{\beta}$ (or \emph{$\beta$-$\alpha$ bordered Heegaard diagram} ${}^{\beta} \mathcal{H}^{\alpha}$); it represents a diagram whose left boundary was parametrized by $\alpha$-arcs and right $\beta$-arcs (left boundary by $\beta$-arcs and right by $\alpha$-arcs). \\

From a doubly bordered Heegaard diagram ${}^{\alpha} \mathcal{H}^{\beta}$, we construct three manifold $Y$ with two boundary components in the following way. Attach two disks $D_L$ and $D_R$ to $\partial_L \overline{\Sigma}$ and $\partial_R \overline{\Sigma}$ to get a closed Heegaard surface $\Sigma$. Then $Y$ is obtained by taking $[0,1] \times \Sigma$ and attach three-dimensional two handles along attaching circles $\boldsymbol{\alpha}^c$ and $\boldsymbol{\beta}^c$. In particular, $[0,1] \times ( D_L \cup ( \textrm{nbd}( \mathbf{z} ) ) \cup D_R )$ is a tunnel connecting two boundaries. \\

\subsection{Half Boundary Dehn Twist}
Next we will discuss a Dehn twist along the boundary. We again repeat the Definitions 3.5 and 3.6 from~\cite{LOT11}.
\begin{defn}
(Definition 3.5, \cite{LOT11}) Let $A$ be an (oriented) annulus with one boundary component marked as the
``inside boundary'' and the other as the ``outside boundary.'' A \emph{radial curve} is any embedded
curve in $A$ which connects the inside and outside boundary of $A$. Suppose that $r$ and $r'$
are two oriented, radial curves which intersect the inside boundary of $A$ at the same point,
but which are otherwise disjoint. We say that $r'$ is \emph{to the right} of $r$ if $r$ has a regular
neighborhood $U$ with an orientation-preserving identification with $(−\epsilon , \epsilon ) \times [0, 1]$, so that $r$ is identified with $ \{ 0 \} \times [0, 1]$, the inside boundary meets $U$ in $(−\epsilon , \epsilon) \times {0}$, and $r' \cap U$ is contained in $[0, \epsilon) \times [0, 1]$.
\end{defn}
\begin{defn}
(Definition 3.6, \cite{LOT11} Let $(F, D_{\alpha} \cup s \cup D_{\beta})$ be a surface with preferred disk decomposed into three parts. A \emph{positive half Dehn twist along the boundary}, denoted $\tau^{1/2}_{\partial}$, is a homeomorphism with the following properties. 
  \begin{itemize}
    \item there is a disk neighborhood $N$ of $D_{\alpha} \cup s \cup D_{\beta}$ so that $\tau^{1/2}_{\partial}$ fixes the complement of $N$.
    \item $\tau^{1/2}_{\partial}$ maps the preferred disk to itself, but switches $D_{\alpha}$ and $D_{\beta}$; and
    \item there is a radial arc $r$ in the annulus $A = N \backslash \textrm{int} (D_{\alpha} \cup s \cup D_{\beta})$ (oriented so that it terminates at $\partial (D_{\alpha} \cup s \cup D_{\beta})$ which is mapped under $\tau^{1/2}_{\partial}$  to a new arc $r'$, which is to the right of $r$. Here, we view $\partial (D_{\alpha} \cup s \cup D_{\beta})$ as the outside boundary of $A$.  
  \end{itemize}   
\end{defn}
Roughly, $\tau^{1/2}_{\partial}$ twists the neighborhood of boundary of $F \backslash D$ $180^{\circ}$ in a ``positive'' direction. Clearly, the composition of two positive half boundary Dehn twists will cause a ``full positive twist'' along the boundary and will be denoted $\tau_{\partial}$. Note that the full boundary Dehn twist acts trivially on a bordered three manifold with single boundary.\\
\begin{figure}
  \centering
  \includegraphics[scale=1]{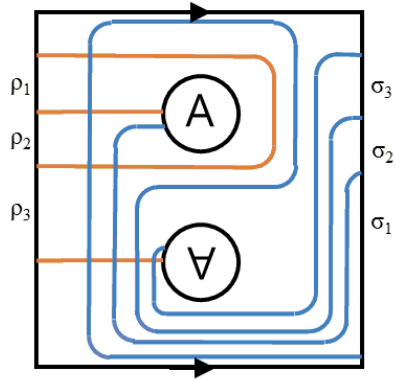} 
  \caption{The genus-1 diagram of ${}^{\alpha} \mathsf{AZ}(T^2)^{\beta}$ interpolating piece. Red lines denote $\alpha$-arcs and blue lines $\beta$-arcs. Gluing two of the diagrams will result in a Heegaard diagram of full boundary Dehn twist.}
  \label{fig:interpolate}
\end{figure}
\subsection{Interpolating piece}
Denis Auroux has first introduced an \emph{interpolating piece} ${}^{\alpha} {\mathsf{AZ}}^{\beta}$ in \cite{Aur10}. For a surface $\mathcal{Z}$ with genus $k$, the interpolating piece, it is an $\alpha$-$\beta$ bordered diagram with genus $k$. For our purpose we will only use the diagram of $k=1$. The diagram can be found in Figure \ref{fig:interpolate}.\\

There are two facts about the interpolating piece, which we will use in this paper. The first proposition is stated below.
\begin{prop}
(\emph{Proposition 4.1,~\cite{LOT11}}) The type $AA$ bimodule $\widehat{CFAA}({}^{\alpha} \mathsf{AZ}( \mathcal{Z})^{\beta} )$ associated with the diagram ${}^{\alpha} \mathsf{AZ} (\mathcal{Z})^{\beta}$, viewed as a left-right $\mathcal{A}( \mathcal{Z} )$-$\mathcal{A}( \mathcal{Z} )$-bimodule, is isomorphic to the bimodule $\mathcal{A}( \mathcal{Z} )$.
\end{prop}
In \cite{LOT10}, type-$AA$ bimodule was defined to be a right-right $\mathcal{A}( \mathcal{Z} )$-$\mathcal{A}( \mathcal{Z} )$-bimodule. However, in this proposition we regard one of its right boundary as its left boundary by taking the opposite algebra $\mathcal{A}( \mathcal{Z} )^{\mathrm{op}} \cong \mathcal{A}( - \mathcal{Z} )$. \\

To properly state the second proposition, we should recall the followings.
\begin{itemize}
  \item (\emph{Constuction 3.21,~\cite{LOT11}.}) Let $Y$ be a strongly bordered three manifold with boundaries parametrized by $\phi_L : F( \mathcal{Z}_L) \rightarrow \partial_L Y$ and $\phi_R : F (\mathcal{Z}_R ) \rightarrow \partial_R Y$. Moreover also suppose $Y$ is homeomorphic to an interval times the surface under a homeomorphism $\Phi : [0,1] \times F( \mathcal{Z}_R ) \rightarrow Y$, such that 1) $\Phi |_{ \{ 1\} \times F (\mathcal{Z}_R) } = \phi_R$, 2) the images under $\Phi$ of $[0,1]$ times the three distinguished regions $D_{\alpha}$, $s$, $D_{\beta}$ (in $F(\mathcal{Z}_R)$) are mapped to the three corresponding distinguished regions in $Y$. Then we can associate $Y$ to a map $\phi_Y : - F( \mathcal{Z}_L ) \rightarrow F( \mathcal{Z}_R)$ by
  \begin{displaymath}
    \phi_Y = ( \Phi |_{ \{ 0 \} \times F (\mathcal{Z}_R ) } )^{-1} \circ (- \phi_L) : F( \mathcal{Z}_L) \rightarrow F( \mathcal{Z}_R).
  \end{displaymath}
  we call $Y$ a \emph{mapping cylinder of $\phi$}.
  \item (\emph{Definition 3.4,~\cite{LOT11}.}) Let $\mathcal{Z}^{\alpha}$ and $\mathcal{Z}^{\beta}$ be two pointed matched circles which only differ in label. In this case, we say that $\mathcal{Z}^{\alpha}$ and $\mathcal{Z}^{\beta}$ are \emph{twin pointed matched circles}. For twin pointed matched circles, there are canonical orientation-reversing homeomorphisms
  \begin{displaymath}
    K_{\alpha, \beta} : F(\mathcal{Z}^{\alpha}) \rightarrow F(\mathcal{Z}^{\beta}) \qquad K_{\beta, \alpha} : F(\mathcal{Z}^{\beta}) \rightarrow F(\mathcal{Z}^{\alpha}).
  \end{displaymath}
  These homeomorphisms preserve $D_{\alpha} \cup D_{\beta} \cup D_{\beta}$, but $D_{\alpha}$ maps to $D_{\beta}$ and vice versa. 
\end{itemize}
\begin{prop}
(\emph{Proposition 4.2,~\cite{LOT11}}) The diagram $\mathsf{AZ}(\mathcal{Z})$ represents a positive half boundary Dehn twist, in the following sense. 
  \begin{displaymath}
    K_{\beta, \alpha} \circ \phi_{ {}^{\alpha} \mathsf{AZ}^{\beta} } = \tau_{\partial}^{1/2} : - F( \mathcal{Z}^{\alpha} ) \rightarrow -F ( \mathcal{Z}^{\alpha} ).
  \end{displaymath}
\label{prop:halftwist}
\end{prop}
We are now ready to prove Theorem 1. \\

\emph{Proof of Theorem 1.} Let ${}^{\mathcal{A'}, \mathcal{A'} } M$ be a (finite dimensional) type-$DD$ module with left-left $dg$ algebra action $\mathcal{A'}= \mathcal{A} ( - \mathcal{Z} )$. Taking the dual of $M$, we get a right-right type-$DD$ module $\overline{M}^{\mathcal{A'}, \mathcal{A'}}$. For any type-$DD$ module ${}^{\mathcal{A'}, \mathcal{A'} } N$, the morphism spaces between them is,
\begin{displaymath}
\overline{M} \boxtimes_{\mathcal{A'}, \mathcal{A'} } (\mathcal{A'} \otimes \mathcal{A'})
\boxtimes_{\mathcal{A'}, \mathcal{A'}} N.
\end{displaymath}
In terms of Heegaard diagrams, Let ${}^{\alpha} \mathcal{H}^{\alpha}_1$ and ${}^{\alpha} \mathcal{H}^{\alpha}_2$ be Heegaard diagrams where $\widehat{CFDD}( \mathcal{H}_1 ) = {}^{\mathcal{A}',\mathcal{A}',}M$ and $\widehat{CFD}( \mathcal{H}_2 ) = {}^{\mathcal{A}',\mathcal{A}'}N$. Then the Heegaard diagram representing the above morphism space is
\begin{displaymath}
\xymatrix{
& {}^{\beta} \overline{ \mathcal{H}_1 } {}^{\beta} & \\
{}^{\alpha} \mathsf{AZ} (\mathcal{Z})^{\beta} \ar@{-}@/^/[ur]|{ {}_{\partial_R} \cup_{\partial_L} } & &
{}^{\beta} \mathsf{AZ} (\mathcal{Z})^{\alpha} \ar@{-}@/_/[ul]|{ {}_{\partial_R} \cup_{\partial_L} } \\
& {}^{\alpha} \mathcal{H}_2^{\alpha}. \ar@{-}@/^/[ul]|{ {}_{\partial_L} \cup_{\partial_L} } \ar@{-}@/_/[ur]|{ {}_{\partial_R} \cup_{\partial_R} } &
}
\end{displaymath}
Now we take advantage of Lemma 4.6,~\cite{LOT11}. According to the lemma, two Heegaard diagrams
\begin{displaymath}
{}^{\alpha} \mathcal{H}^{\mathcal{\alpha}} {}_{\partial_R} \cup_{\partial_L} {}^{\alpha} \mathsf{AZ} (- \mathcal{Z})^{\beta}
\end{displaymath}
and
\begin{displaymath}
{}^{\alpha} \mathsf{AZ} (\mathcal{Z})^{\beta} {}_{\partial_R} \cup_{\partial_L} {}^{\beta} (- \overline{ \mathcal{H} } )^{\mathcal{\beta}}
\end{displaymath}
represent the same strongly bordered three-manifold ( ${}^{\beta} (- \overline{ \mathcal{H} } )^{\mathcal{\beta}}$ means a Heegaard diagram obtained from $\mathcal{H}$, by switching $\alpha$ and $\beta$ labels and choosing the opposite orientation from $\mathcal{H}$. Note that $Y(-\overline{H})$ is an orientation-preserving homeomorphic to $Y(\mathcal{H})$ ). Then the above diagram is
\begin{displaymath}
\xymatrix{
& {}^{\beta} \overline{ \mathcal{H}_1 } {}^{\beta} & \\
{}^{\alpha} \mathsf{AZ} (\mathcal{Z})^{\beta} \ar@{-}@/^/[ur]|{ {}_{\partial_R} \cup_{\partial_L} } & &
{}^{\beta} (- \overline{\mathcal{H}_2} )^{\beta} \ar@{-}@/_/[ul]|{ {}_{\partial_R} \cup_{\partial_L} } \\
& {}^{\alpha} \mathsf{AZ}( - \mathcal{Z})^{\beta}. \ar@{-}@/^/[ul]|{ {}_{\partial_L} \cup_{\partial_L} } \ar@{-}@/_/[ur]|{ {}_{\partial_R} \cup_{\partial_R} } &
}
\end{displaymath}
By Proposition~\ref{prop:halftwist}, an interpolating piece represents a half boundary Dehn Twist. Thus, two consecutive interpolating pieces represent a negative full boundary Dehn twist $\tau^{-1}_{\partial} : F ( \mathcal{Z} ) \rightarrow F(\mathcal{Z}) $ (See Corollary 4.5,~\cite{LOT11}). This proves the morphism space $\mathrm{Mor}^{\mathcal{A'}, \mathcal{A'}} (\widehat{CFDD}(Y_1), \widehat{CFDD}(Y_2))$ is isomorphic to the $\widehat{CF}$ of the above diagram, or equivalently
\begin{displaymath}
\widehat{CFAA} ( \tau^{-1}_{\partial} ( Y(\mathcal{H}_2) ) \boxtimes_{\mathcal{A}, \mathcal{A}'} \widehat{CFDD} ( -Y(\mathcal{H}_1) ). \ \ \square
\end{displaymath}

\begin{figure}
  \centering
  \includegraphics[scale=1]{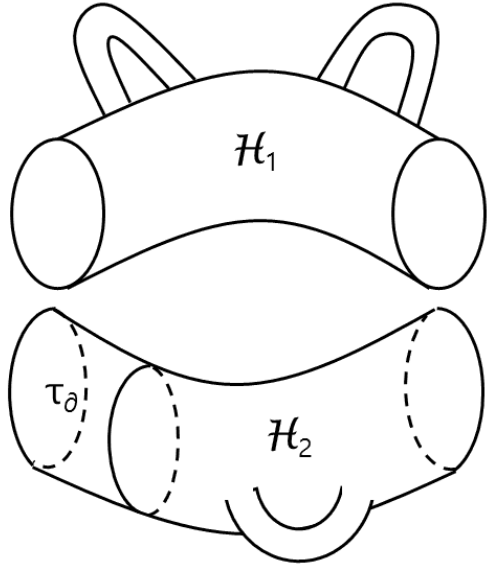}
  \caption{The Heegaard Floer homology $\widehat{CF}$ of above Heegaard diagram will result in a chain complex homotopic to the morphism space. The $\alpha$- $\beta$ labelling of diagram $\mathcal{H}_1$ has to be interchanged, and for $\mathcal{H}_2$, the orientation of the diagram should be reversed, as well as the labelling.} 
\end{figure}

\section{Morphism space $\textrm{Mor}(\widehat{CFDD}(\mathbb{I}), \widehat{CFDD}({\tau_{m+1}}) )$}
We consider type-$DD$ morphisms with two left actions of torus algebras $\mathcal{A}(-T^2), \mathcal{A} (-T^2)$. We will follow the common convention; $\rho$'s and $\sigma$'s will denote ``left'' and ``right'' algebra elements, respectively. \\
From now on, we will denote $\mathcal{A} = \mathcal{A}(T^2)$ and $\mathcal{A}' = \mathcal{A}(-T^2)$, unless otherwise specified.\\

First, we recall the two type-$DD$ structures of our concern. The identity type-$DD$ module $\widehat{CFDD} (\mathbb{I})$ has two generators $\mathbf{x}$ and $\mathbf{y}$ with idempotent actions $\iota_0 \jmath_0 \mathbf{x} = \mathbf{x}$, $\iota_1 \jmath_1 \mathbf{y} = \mathbf{y}$, and otherwise zero. \\

The differential $(\delta^I)^1$ was shown to be
\begin{eqnarray*}
  (\delta^I)^1 \mathbf{x} & = & (\rho_1 \sigma_3 + \rho_3 \sigma_1 + \rho_{123} \sigma_{123} ) \mathbf{y} \\
  (\delta^I)^1 \mathbf{y} & = & \rho_2 \sigma_2 \mathbf{x}.
\end{eqnarray*}

The computation of $\widehat{CFDD} ( \tau_{m+1} )$ can be easily done, by dualizing one end of $\widehat{CFDA} ( \tau_{m+1} )$, which was computed in~\cite{LOT08}. It has three generators $\mathbf{p}, \mathbf{q}$ and $\mathbf{r}$ such that
\begin{eqnarray*}
\iota_0 \jmath_0 \mathbf{p} & = & \mathbf{p} \\
\iota_1 \jmath_1 \mathbf{q} & = & \mathbf{q} \\
\iota_1 \jmath_0 \mathbf{r} & = & \mathbf{r},
\end{eqnarray*}
and otherwise zero. \\

$\widehat{CFDA} ( \tau_{m+1} )$ has the following nine nontrivial differentials listed as below.
\begin{eqnarray*}
(\mathbf{p}, \rho_1) & \mapsto & \rho_1 \mathbf{q} \\
(\mathbf{p}, \rho_{123}) & \mapsto & \rho_{123} \mathbf{q} \\
(\mathbf{p}, \rho_3, \rho_{23}) & \mapsto & \rho_3 \mathbf{q} \\
(\mathbf{p}, \rho_{12}) & \mapsto & \rho_{123} \mathbf{r} \\
(\mathbf{p}, \rho_3, \rho_2) & \mapsto & \rho_3 \mathbf{r} \\
(\mathbf{q}, \rho_{23}) & \mapsto & \rho_{23} \mathbf{q} \\
(\mathbf{q}, \rho_2) & \mapsto & \rho_{23} \mathbf{q} \\
(\mathbf{r}, \rho_3) & \mapsto & \mathbf{q} \\
(\mathbf{r}, 1) & \mapsto & \rho_2 \mathbf{q}.
\end{eqnarray*}
Dualizing the third differential will cause algebra elements whose endpoints do not match, so it will be zero. In addition, dualizing the fourth and sixth differential will be zero due to idempotents. Thus the differential $(\delta^{\tau_{m+1}})^1$ is dualized to
\begin{eqnarray*}
(\delta^{\tau_{m+1}})^1 (\mathbf{p}) & = & \rho_1 \sigma_3 \mathbf{q} + \rho_{123} \sigma_{123} \mathbf{q} + \rho_3 \sigma_{12} \mathbf{r} \\
(\delta^{\tau_{m+1}})^1 (\mathbf{q}) & = & \rho_{23} \sigma_2 \mathbf{r} \\
(\delta^{\tau_{m+1}})^1 (\mathbf{r}) & = & \rho_2 \mathbf{p} + \sigma_1 \mathbf{q}.
\end{eqnarray*}

Since every morphism in $\textrm{Mor}^{\mathcal{A}(-T^2), \mathcal{A}(-T^2)} (\widehat{CFDD}(\mathbb{I}), \widehat{CFDD}({\tau_{m+1}}) )$ has to respect the idempotent actions, $\mathbf{x}$ and $\mathbf{y}$ should be mapped to 
\begin{displaymath}
\mathbf{x} \mapsto \left\{
\begin{array}{ccc}
  a_1 \cdot \rho_1 \mathbf{r} & +\ a_7 \cdot \mathbf{p} & +\ a_{11} \cdot \rho_1 \sigma_1 \mathbf{q} \\ 
  +\ a_2 \cdot \rho_3 \mathbf{r} & +\ a_8 \cdot \rho_{12} \mathbf{p} & +\ a_{12} \cdot \rho_1 \sigma_3 \mathbf{q} \\ 
  +\ a_3 \cdot \rho_1 \sigma_{12} \mathbf{r} & +\ a_9 \cdot \sigma_{12} \mathbf{p} & +\ a_{13} \cdot \rho_3 \sigma_1 \mathbf{q} \\ 
  +\ a_4 \cdot \rho_3 \sigma_{12} \mathbf{r} & +\ a_{10} \cdot \rho_{12} \sigma_{12} \mathbf{p} & +\ a_{14} \cdot \rho_3 \sigma_3 \mathbf{q} \\ 
  +\ a_5 \cdot \rho_{123} \mathbf{r} & & +\ a_{15} \cdot \rho_{123} \sigma_1 \mathbf{q} \\ 
  +\ a_6 \cdot \rho_{23} \sigma_{12} \mathbf{r} & & +\ a_{16} \cdot \rho_{123} \sigma_3 \mathbf{q} \\ 
   & & +\ a_{17} \cdot \rho_1 \sigma_{123} \mathbf{q} \\ 
   & & +\ a_{18} \cdot \rho_3 \sigma_{123} \mathbf{q} \\ 
   & & +\ a_{19} \cdot \rho_{123} \sigma_{123} \mathbf{q},
\end{array}
\right.
\end{displaymath}

\begin{displaymath}
\mathbf{y} \mapsto \left\{
\begin{array}{ccc}
  b_1 \cdot \sigma_2 \mathbf{r} & +\ b_3 \cdot \mathbf{q} & +\ b_7 \cdot \rho_2 \sigma_2 \mathbf{p} \\
  +\ b_2 \cdot \rho_{23} \sigma_2 \mathbf{r} & +\ b_4 \cdot \sigma_{23} \mathbf{q} & \\
  & +\ b_5 \cdot \rho_{23} \mathbf{q} & \\
  & +\ b_6 \cdot \rho_{23} \sigma_{23} \mathbf{q}. &  
\end{array}
\right.
\end{displaymath}
where $a_i$s and $b_j$s are $\mathbb{F}_2$ coefficients. \\

The morphism space $\textrm{Mor}_{\mathcal{A}', \mathcal{A}'} (\widehat{CFDD}(\mathbb{I}), \widehat{CFDD}({\tau_{m+1}}) )$ has an obvious differential, and for $f$ to be in the kernel of it, $f$ must be a chain map. In other words, $(\delta^{\tau_{m+1}})^1 f \mathbf{x} = f (\delta^I)^1 \mathbf{x}$ and $(\delta^{\tau_{m+1}})^1 f \mathbf{y} = f (\delta^I)^1 \mathbf{y}$. Then we get the following equations.

\begin{eqnarray*}
(\delta^{\tau_{m+1}})^1 f \mathbf{x} & = & (\delta^{\tau_{m+1}})^1 ( a_1 \cdot \rho_1 \mathbf{r} + \cdots + a_{19} \cdot \rho_{123} \sigma_{123} \mathbf{q} ) \\
& = & ( a_1 \cdot \rho_1 + \cdots + a_6 \cdot \rho_{123} \sigma_{12} ) \cdot ( \sigma_1 \mathbf{q} + \rho_2 \mathbf{p} ) \\
& & + ( a_7 + a_8 \cdot \rho_{12} + \cdots + a_{10} \cdot \rho_{12} \sigma_{12} ) \cdot ( \rho_1 \sigma_3 \mathbf{q} + \rho_{123} \sigma_{123} \mathbf{q} + \rho_3 \sigma_{12} \mathbf{r} ) \\
& & + ( a_{11} \cdot \rho_1 \sigma_1 + \cdots + a_{19} \cdot \rho_{123} \sigma_{123} ) \cdot ( \rho_{23} \sigma_2 \mathbf{r} ) \\
& = & a_1 \cdot ( \rho_1 \sigma_1 \mathbf{q} + \rho_{12} \mathbf{p} ) + a_2 \cdot \rho_3 \sigma_1 \mathbf{q} + a_3 \cdot \rho_{12} \sigma_{12} \mathbf{p} + a_5 \cdot \rho_{123} \sigma_1 \mathbf{q} \\
& & + a_7 \cdot ( \rho_1 \sigma_3 \mathbf{q} + \rho_{123} \sigma_{123} \mathbf{q} + \rho_3 \sigma_{12} \mathbf{r} ) + a_8 \cdot \rho_{123} \sigma_{12} \mathbf{r} + a_9 \cdot \rho_1 \sigma_{123} \mathbf{q} \\
& & + a_{11} \cdot \rho_{123} \sigma_{12} \mathbf{r}, 
\end{eqnarray*}
\begin{eqnarray*}
f (\delta^I)^1 \mathbf{x} & = & f ( \rho_{123} \sigma_{123} \mathbf{y} + \rho_1 \sigma_3 \mathbf{y} + \rho_3 \sigma_1 \mathbf{y} ) \\
& = & ( \rho_{123} \sigma_{123} + \rho_1 \sigma_3 + \rho_3 \sigma_1 ) \\
& & \cdot ( (b_1 \cdot \sigma_2 + b_2 \cdot \rho_{23} \sigma_2) \mathbf{r} + (b_3 + b_4 \cdot \sigma_{23} + b_5 \cdot \rho_{23} + b_6 \cdot \rho_{23} \sigma_{23} ) \mathbf{q} + b_7 \cdot \rho_2 \sigma_2 \mathbf{p} ) \\
& = & b_1 \cdot \rho_3 \sigma_{12} \mathbf{r} + b_3 \cdot ( \rho_{123} \sigma_{123} + \rho_1 \sigma_3 + \rho_3 \sigma_1 ) \mathbf{q} + b_4 \cdot \rho_3 \sigma_{123} \mathbf{q} \\
& & + b_5 \cdot \rho_{123} \sigma_3 \mathbf{q}. 
\end{eqnarray*}
These equations force $a_1 = a_3 = a_5 = a_9 =0$, $a_8 = a_{11}$, $b_4 = b_5 =0$, and $a_2 = a_7 = b_1 = b_3$. Similarly,
\begin{eqnarray*}
(\delta^{\tau_{m+1}})^1 f \mathbf{y} & = & (\delta^{\tau_{m+1}})^1 ( b_1 \cdot \sigma_2 \mathbf{r} + b_2 \cdot  \rho_{23} \sigma_2 + b_3 \cdot \mathbf{r} + b_6 \cdot \mathbf{q} + b_7 \cdot \rho_{23} \sigma_{23} ) \mathbf{q} \\
& = & b_1 \cdot \rho_2 \sigma_2 \mathbf{p} + b_3 \cdot \rho_{23} \sigma_2 \mathbf{r},
\end{eqnarray*}
\begin{eqnarray*}
f (\delta^I)^1 \mathbf{y} & = & \rho_2 \sigma_2 f(\mathbf{x}) \\
& = & \rho_2 \sigma_2 ( a_2 \cdot \rho_3 \mathbf{r} + \cdots + a_{19} \cdot \rho_{123} \sigma_{123} \mathbf{q} ) \\
& = & a_2 \cdot \rho_{23} \sigma_2 \mathbf{r} + a_7 \cdot \rho_2 \sigma_2 \mathbf{p} + a_{14} \cdot \rho_{23} \sigma_{23} \mathbf{q}. 
\end{eqnarray*}
Thus $a_{14} =0$. \\

It is clear that $\textrm{Mor}_{\mathcal{A}', \mathcal{A}'} (\widehat{CFDD}(\mathbb{I}), \widehat{CFDD}({\tau_{m+1}}) )$ is a 26-dimensional vector space, and we can describe any type-$DD$ morphism  
by determining the coefficients $a_i$ or $b_j$. From now on, we will denote the map as $\sum (a_i) + \sum (b_j)$ without loss of information.\\ 

We are also interested in the image of differential $\partial$ of the morphism space $\textrm{Mor}_{\mathcal{A}', \mathcal{A}'} (\widehat{CFDD}(\mathbb{I}), \widehat{CFDD}({\tau_{m+1}}) )$, i.e, maps homotopic to zero. For any morphism $H \in \textrm{Mor}_{\mathcal{A}', \mathcal{A}'} (\widehat{CFDD}(\mathbb{I}), \widehat{CFDD}({\tau_{m+1}}) )$, the differential $\partial$ is defined to be 
\begin{displaymath}
\partial H = (\delta^{\tau_{m+1}})^1 \circ H + H \circ (\delta^{I})^1.
\end{displaymath} 
Let $H = (a_i)$ or $(b_j)$ for some $i$ and $j$. By a straightforward computation,
\begin{eqnarray*}
\partial (a_1) & = & (a_8) + (a_{11}) \\
\partial (a_2) & = & (a_{13}) + (b_2) \\
\partial (a_3) & = & (a_{10}) \\
\partial (a_5) & = & (a_{15}) \\
\partial (a_7) & = & (a_4) + (a_{12}) + (a_{19}) + (b_7) \\
\partial (a_8) & = & \partial (a_{11}) = (a_6) \\
\partial (a_9) & = & (a_{17}) \\
\partial (a_{14}) & = & (b_6) \\
\partial (b_1) & = & (a_4) + (b_7) \\
\partial (b_3) & = & (a_{12}) + (a_{13}) + (a_{19}) + (b_2) \\
\partial (b_4) & = & (a_{18}) \\
\partial (b_5) & = & (a_{16}),  
\end{eqnarray*}
otherwise $\partial H = 0$. \\

Now we can explicitly describe the basis of $H_* ( \textrm{Mor}_{\mathcal{A}', \mathcal{A}'} (\widehat{CFDD}(\mathbb{I}), \widehat{CFDD}({\tau_{m+1}}) ) )$; it is a four dimensional space generated by $(a_2) + (a_7) + (b_1) + (b_3)$, $(a_4)$, $(a_{12})$, and $(a_{13})$.  

\section{From $\widehat{CFDD}$ map to $\widehat{CFDA}$ map}
In this section, we compute type-$DA$ morphisms from $\widehat{CFDD}(\mathbb{I})$ to $\widehat{CFDD}({\tau_{m+1}}) ) )$. \\

In the previous section, we explicitly computed type-$DD$ morphism spaces from $\widehat{CFDD}(\mathbb{I})$ to $\widehat{CFDD}({\tau_{m+1}}) ) )$. In~\cite{LOT10}, the functor $\cdot \boxtimes {}_{\mathcal{A}} [ \mathbb{I}]^{\mathcal{A}} : {}_{\mathcal{A}, \mathcal{A}} \mathsf{Mod} \rightarrow {}_{\mathcal{A}} \mathsf{Mod}_{\mathcal{A}}$ from the category of left-left type-$DD$ bimodule to the category of left-right type-$DA$ bimodule, is proven to be categorical equivalence. By taking the $\mathcal{A}_{\infty}$ tensor product with the identity bimodule ${}^{\mathcal{A}} \mathbb{I}_{\mathcal{A}}$, the type-$DD$ morphisms can be transformed into type-$DA$ morphisms. For the model of the identity bimodule, we will use a Heegaard diagram $\mathbb{I}'$ that was given in the section 10 of \cite{LOT10} to compute the identity $DA$ module. Of course we may use the Heegaard diagram $\mathbb{I}$ with only two generators, but it is not operationally bounded in the sense of \cite{LOT10} (i.e, it has infinitely many nontrivial differentials) so there are infinitely many relations to consider. \\

First, we should recall the type-$AA$ identity bimodule $[\mathbb{I}']_{\mathcal{A}, \mathcal{A}}$ introduced in the section 10 of \cite{LOT10}. The bimodule has five generators $\mathbf{x}$, $\mathbf{y}$, $\mathbf{w_1}$, $\mathbf{w_2}$, $\mathbf{z_1}$ and $\mathbf{z_2}$, with the following $\mathcal{A}_{\infty}$ relations.
\begin{itemize}
\item $(\mathbf{w_1}, \sigma_1) = \mathbf{y}$ 
\item $(\mathbf{w_1}, \sigma_{12}, \rho_2) = \mathbf{x}$
\item $(\mathbf{w_1}) = (\mathbf{w_1}, \sigma_{12}, \rho_{23}) = \mathbf{w_2}$
\item $(\mathbf{w_1}, \sigma_{123}, \rho_2) = (\mathbf{w_2}, \sigma_3, \sigma_2, \sigma_1, \rho_2) = \mathbf{z_2}$
\item $(\mathbf{z_1}, \rho_1) = \mathbf{y}$ 
\item $(\mathbf{z_1}, \sigma_2, \rho_{12}) = \mathbf{x}$
\item $(\mathbf{z_1}) = (\mathbf{z_1}, \sigma_{23}, \rho_{12}) = \mathbf{z_2}$
\item $(\mathbf{z_1}, \sigma_2, \rho_{123}) = \mathbf{w_2}$
\item $(\mathbf{y}, \sigma_2, \rho_2) = \mathbf{x}$
\item $(\mathbf{y}, \sigma_2, \rho_{23}) = \mathbf{w_2}$
\item $(\mathbf{y}, \sigma_{23}, \rho_2) = \mathbf{z_2}$
\item $(\mathbf{x}, \rho_3) = \mathbf{w_2}$
\item $(\mathbf{x}, \sigma_3) = \mathbf{z_2}$
\end{itemize}

To the left hand side of the above type-$AA$ bimodule, we take a tensor product with type-$DD$ identity bimodule ${}^{\mathcal{A}, \mathcal{A}} [\mathbb{I}]$ with two generators $\mathbf{x}$ and $\mathbf{y}$, with $\delta^1 \mathbf{x} = (\rho_1 \sigma_3 + \rho_3 \sigma_1 + \rho_{123} \sigma_{123} ) \mathbf{y}$ and  $\delta^1 \mathbf{y} = \rho_2 \sigma_2 \mathbf{x}$. Then, the type-$DA$ bimodule $( {}^{\mathcal{A}, \mathcal{A}} [\mathbb{I}]) \boxtimes ([\mathbb{I}']_{\mathcal{A}, \mathcal{A}} )$ has the following $DA$ structure.

\begin{displaymath}
\xymatrix@C=2cm@R=1.5cm{
\mathbf{xw_1} \ar[dr]|{\rho_3 \otimes 1} \ar[ddd]_1 & & \mathbf{yz_1} \ar[dl]|{1 \otimes \rho_1} \ar[ddd]^1 \ar@/^/[ddl]|{\rho_2 \otimes \rho_{12} } \ar@/^5cm/[dddll]|{\rho_2 \otimes \rho_{123}}\\
& \mathbf{yy} \ar[d]|{\rho_2 \otimes \rho_2} \ar@/_/[ddl]|{\rho_2 \otimes \rho_{23}} & \\
& \mathbf{xx} \ar[dl]|{1 \otimes \rho_3} \ar[dr]|{\rho_1 \otimes 1} & \\
\mathbf{xw_2} & & \mathbf{y_2} 
}
\end{displaymath} 
\bigskip

For computational convenience, we take the $\mathcal{A}_{\infty}$ tensor product of $\widehat{CFDD} ( \tau_{m+1} )$ with $[\mathbb{I}']_{\mathcal{A}, \mathcal{A}}$. The result is given below.

\begin{displaymath}
\xymatrix@C=2cm@R=1.5cm{
& \mathbf{qy} \ar[d]|{\rho_{23} \otimes \rho_{23}} \ar[dr]|{\rho_{23} \otimes \rho_2} & \mathbf{qz_1} \ar[r]|1 \ar[d]|{\rho_{23} \otimes \rho_{12}} \ar[l]|{1 \otimes \rho_1} \ar@/_2.7cm/[dl]|{\rho_{23} \otimes \rho_{123}} & \mathbf{qz_2} \\
\mathbf{rw_1} \ar[ur]|{1} \ar[r]|{1} \ar[d]|(0.3){\rho_2 \otimes 1}& \mathbf{rw_2} \ar[d]|(0.3){\rho_2 \otimes 1} & \mathbf{rx} \ar[l]|{1 \otimes \rho_3} \ar[d]|(0.3){\rho_2 \otimes 1} & \\
\mathbf{pw_1} \ar[r]|{1} \ar[ur]|{\rho_3 \otimes \rho_{23}} \ar[urr]|(0.7){\rho_3 \otimes \rho_2 } & \mathbf{pw_2} & \mathbf{px} \ar[l]|{1 \otimes \rho_3} \ar[uur]|{\rho_1 \otimes 1} & \\
}
\end{displaymath}
Finally, we compute homologically nontrivial type-$DA$ maps from $( {}^{\mathcal{A}, \mathcal{A}} [\mathbb{I}]) \boxtimes ([\mathbb{I}']_{\mathcal{A}, \mathcal{A}} )$ to $\widehat{CFDD}(\tau_{m+1}) \boxtimes [\mathbb{I}']_{\mathcal{A}, \mathcal{A}}$. Four nontrivial type-$DD$ morphisms are converted as below. 
\begin{itemize}
  \item $(a_2)+(a_7)+(b_1)+(b_3)$ is interpreted as $f_1( \mathbf{x} ) = \mathbf{p} + \rho_3 \mathbf{r}, \ f_1( \mathbf{y} ) = \mathbf{q} + \sigma_2 \mathbf{r}$. Then $(f_1 \boxtimes \mathbb{I}')$ is
    \begin{eqnarray*}
      (\mathbf{yy}, \rho_2) & = & \mathbf{rx} \\
      (\mathbf{yz_1}, \rho_{12}) & = & \mathbf{rx} \\
      (\mathbf{yy}, \rho_{23}) & = & \mathbf{rw_2} \\
      (\mathbf{yz_1}, \rho_{123}) & = & \mathbf{rw_2} \\
      (\mathbf{xw_1}, \rho_2) & = & \rho_1 \mathbf{qz_2} \\
      (\mathbf{xw_1} ) & = & \mathbf{pw_1} + \rho_3 \mathbf{rw_1} \\
      (\mathbf{xw_2} ) & = & \mathbf{pw_2} + \rho_3 \mathbf{rw_2} \\
      (\mathbf{xx} ) & = & \mathbf{px} + \rho_3 \mathbf{rx} \\
      (\mathbf{yy} ) & = & \mathbf{qy} \\
      (\mathbf{yz_1} ) & = & \mathbf{qz_1} \\
      (\mathbf{yz_2} ) & = & \mathbf{qz_2}.
    \end{eqnarray*}
    \item $(a_4)$ is $f_2 ( \mathbf{x} ) = \rho_3 \sigma_{12} \mathbf{r}, \ f_2 ( \mathbf{y} ) =0 $. Then $(f_2 \boxtimes \mathbb{I}')$ is
    \begin{eqnarray*}
       (\mathbf{xw_1}, \rho_2) & = & \rho_3 \mathbf{rx} \\
       (\mathbf{xw_1}, \rho_{23}) & = & \rho_3 \mathbf{rw_2}.
    \end{eqnarray*}
    \item $(a_{12})$ is $f_3 ( \mathbf{x} ) = \rho_1 \sigma_3 \mathbf{q}, \ f_3 ( \mathbf{y} ) =0 $. Then $(f_3 \boxtimes \mathbb{I}')$ is
    \begin{eqnarray*}
       (\mathbf{xx}) & = & \rho_1 \mathbf{qz_2} \\
       (\mathbf{xw_1}, \rho_2) & = & \rho_{123} \mathbf{qz_2}.
    \end{eqnarray*}
    \item $(a_{13})$ is $f_4 ( \mathbf{x} ) = \rho_3 \sigma_1 \mathbf{q}, \ f_4 ( \mathbf{y} ) =0 $. Then $(f_4 \boxtimes \mathbb{I}')$ is
    \begin{eqnarray*}
       (\mathbf{xw_1}) & = & \rho_3 \mathbf{qy} \\
       (\mathbf{xw_1}, \rho_2) & = & \rho_{123} \mathbf{qz_2}.
    \end{eqnarray*}
\end{itemize}
For notational simplicity, we will denote $(f_i \boxtimes \mathbb{I}') =: f_i$ for the following sections.
\begin{remark}
Our discussion regarding the morphism space was purely algebraic. For a geometric aspect, an interested reader should refer to Tova Brown's thesis~\cite{Br11}, which includes the discussion on the bordered Floer homology and the mapping cylinder of a map between two Riemann surfaces associated with pointed matched circles. We may relate this discussion to a Lefschetz fibration of a square with a singular point, whose fiber is $T^2$.
\end{remark}

\section{Type-$D$ structure on Direct Limit}
From now on, we will assume the knot discussed in this paper is equipped with a (sufficiently large) positive framing.\\

The type-$D$ module associated with a knot complement $S^3 \backslash \textrm{nbd}(K)$ is well studied in the Chapter 11 of \cite{LOT08}. We will review the basics in bordered Floer homologies of torus boundary case. \\

As a vector space, $\widehat{CFD}(S^3 \backslash \textrm{nbd}(K))$ consists of two summands depending on idempotents $\iota_0$ and $\iota_1$, and as a vector space, it has the following decomposition.
\begin{displaymath}
  \widehat{CFD}(S^3 \backslash \textrm{nbd}(K)) = \iota_0 \widehat{CFD}(S^3 \backslash \textrm{nbd}(K)) \oplus \iota_1 \widehat{CFD}(S^3 \backslash \textrm{nbd}(K)).
\end{displaymath}
The differential $\delta^1$ of $\widehat{CFD}(S^3 \backslash \textrm{nbd}(K))$ can be decomposed as well, depending on torus algebra elements. More precisely, we temporarily view $\widehat{CFD}(S^3 \backslash \textrm{nbd}(K))$ as a vector space and forget the algebra element of differential $\delta^1$. Then the differential called the \emph{coefficient map} $D_I$, where $I=\{ \varnothing, 1, 2, 3, 12, 23, 123 \}$, and the original differential $\delta^1$ are recovered as
\begin{displaymath}
  \delta^1 = 1 \otimes D_{\varnothing} + \rho_1 \otimes D_1 + \rho_2 \otimes D_2 + \rho_3 \otimes D_3 + \rho_{12} \otimes D_{12} + \rho_{23} \otimes D_{23} + \rho_{123} \otimes D_{123}.
\end{displaymath}
For an appropriate choice of the labelling(11.4, \cite{LOT08}) of Heegaard diagram of $S^3 \backslash \textrm{nbd}(K)$, it is well known that $(\iota_0 \widehat{CFD}(S^3 \backslash \textrm{nbd}(K)), D_{\varnothing} )$ is homotopy equivalent to the knot Floer homology $\widehat{CFK}(K)$. \\

To properly state the main result of the Chapter 11 of \cite{LOT08}, we will need the following definitions as well. 
\begin{defn}(Definition 11.25, \cite{LOT08})
A knot Floer homology $CFK^- (K)$ is \emph{reduced} if its differential $\partial$ strictly decreases either the Alexander filtration or $U$-filtration.  
\label{def:reduced}
\end{defn}

In the knot Floer homology, the most commonly used convention is to put generators of the knot Floer complex on the plane with integral coordinate, so that the difference of $x$ coordinate represents difference of $U$ filtration level and $y$ coordinate represents the Alexander filtration. Under the convention the arrow represents a nontrivial differential. Since $CFK^-(K)$ has Alexander and $U$-filtrations, a knot Floer complex is reduced implying there is no nontrivial differential staying at the same filtration level. In other words, every arrow goes to the left or bottom, or diagonally to the bottom-left direction. \\
\begin{defn}(Definition 11.24, \cite{LOT08})
Let $\cdots \subset V_i \subset V_{i+1} \subset \cdots$ be an ascending filtration of a vector space $V$, so $V = \bigcup^{\infty}_{i= - \infty} V_i$. Assume that $\bigcap^{\infty}_{i=-\infty} V_i =0$. Given $v \in V$, define $A(v) := \textrm{inf} \{ i \in \mathbb{Z} | v \in V_i \}$; we call $A(v)$ the \emph{filtration level of $v$}. Let $\textrm{gr}(V)$ denote the associated graded vector space, i.e.,
\begin{displaymath}
\textrm{gr}(V) := \bigoplus_{i=-\infty}^{\infty} \textrm{gr}_i(V) \qquad \textrm{gr}_i(V) := V_i / V_{i-1}.
\end{displaymath} 
Composing the projection $V_i \rightarrow V_i / V_{i-1}$ and the inclusion $V_i / V_{i-1} \rightarrow \textrm{gr}(V)$ gives a map $ \iota_i : V_i \rightarrow \textrm{gr}(V)$. For $v \in V$, define $[v] := \iota_{A(v)}(v)$. \\
A \emph{filtered basis} of a filtered vector space $V$ is a basis $\{ v_1, \cdots v_n \}$ for $V$ such that $\{ [v_1], \cdots, [v_n] \}$ is a basis for $\textrm{gr}(V)$.
\end{defn}
To take computational advantage, we will need to choose special generators of the knot Floer complex. A filtered basis $\{ \xi_i \}$ over $\mathbb{F}[U]$ for $C$ is \emph{vertically simplified} if for each basis vector $\xi_i$, either $\partial \xi_i \in U \cdot C$ or $\partial \xi_i \equiv \xi_{i+1}$ (mod $U \cdot C$). In the latter case, the Alexander filtration difference between $\xi_i$ and $\xi_{i+1}$ is called the \emph{length of the arrow}. A filtered basis $\{ \eta_i \}$ over $\mathbb{F}[U]$ for $C$ is \emph{horizontally simplified} if for each basis vector $\eta_i$, either $A(\partial \eta_i) < A(\eta_i)$ or there is an $m$ so that $\partial \eta_i = U^m \cdot \eta_{i+1} + \epsilon$ where $A(\eta_i) = A(U^m \cdot \eta_{i+1})$ and $A(\epsilon) < A(\eta_i)$. Again, in the latter case, the $U$-filtration difference between $\eta_i$ and $\eta_{i+1}$ is called the \emph{length of the arrow}, too.\\

Now we can spell out the main result of the Chapter 11 of \cite{LOT08}. Note that every finitely generated complex is homotopy equivalent to the reduced complex, and we can always choose a horizontally and vertically simplified basis. Let us suppose $\iota_0 \widehat{CFD}(S^3 \backslash \textrm{nbd}(K)) (= \widehat{CFK}(K))$ is reduced. For any pair $\xi_i, \xi_{i+1} \in \iota_0 \widehat{CFD}(S^3 \backslash \textrm{nbd}(K))$ whose (vertical)length of the arrow is $l$, there is a sequence of elements $\kappa^i_1, \cdots \kappa^i_l \in \iota_1 \widehat{CFD}(S^3 \backslash \textrm{nbd}(K))$ (up to homotopy)such that,
\begin{displaymath}
\xymatrix{
\xi_i \ar[r]^{\rho_1} & \kappa^i_1 & \cdots \ar[l]_{\rho_{23}} & \kappa^i_k \ar[l]_{\rho_{23}} & \kappa^i_{k+1} \ar[l]_{\rho_{23}} & \cdots \ar[l]_{\rho_{23}} & \kappa^i_l \ar[l]_{\rho_{23}} & \xi_{i+1} \ar[l]_{\rho_{123}}.
}
\end{displaymath}
Similarly, for any pair $\eta_i, \eta_{i+1} \in \iota_0 \widehat{CFD}(S^3 \backslash \textrm{nbd}(K))$ whose (horizontal)length of the arrow is $l$, there is a sequence of elements $\lambda^i_1, \cdots, \lambda^i_l \in \iota_1 \widehat{CFD}(S^3 \backslash \textrm{nbd}(K))$(up to homotopy) such that
\begin{displaymath}
\xymatrix{
\eta_i \ar[r]^{\rho_3} & \lambda^i_1 \ar[r]^{\rho_{23}} & \cdots \ar[r]^{\rho_{23}} & \lambda^i_k \ar[r]^{\rho_{23}} & \lambda^i_{k+1} \ar[r]^{\rho_{23}} & \cdots \ar[r]^{\rho_{23}} & \lambda^i_l \ar[r]^{\rho_2} & \eta_{i+1}.
}
\end{displaymath}
For a vertically simplified basis $\{ \xi_i \}$, there is a distinguished element in $\{ \xi_i \}$, denoted $\xi_0$ and is defined as follows. A \emph{vertical complex} $C^v : = C / U \cdot C$, inherits the Alexander filtration, has its homology $H^v (C)$. Since we are considering the knot in $S^3$, the homology has rank 1. Because we are assuming that our basis of $C$ is simplified, there is an element $\xi_0 \in \{ \xi_i \}$ which is the generator of $H^v(C)$. Likewise, we can choose a distinguished element $\eta_0 \in \{ \eta_i \}$ by considering the horizontal complex. Then there exists a sequence of elements $\{ \mu_i \}$ in $\iota_1 \widehat{CFD}(S^3 \backslash \textrm{nbd}(K))$ such that
\begin{displaymath}
\xymatrix{
\xi_0 \ar[r]^{\rho_{123}} & \mu_1 \ar[r]^{\rho_{23}} & \cdots \ar[r]^{\rho_{23}} & \mu_m \ar[r]^{\rho_2} & \eta_0.
}
\end{displaymath} 
\begin{remark}
The length $m$ of the sequence $\{ \mu_i \}$ is determined by the concordance invariant $\tau$ for knots. The invariant $\tau(K)$ is defined to be the minimal $s$ for which the generator of $H_* (\widehat{CF}(S^3))$ can be represented as a sum of generators in Alexander grading less than or equal to $s$. In fact, $m=n-2\tau(K)$, where $n$ is the framing of the knot. \\
\end{remark}
We now associate a morphism $f \in \textrm{Mor}_{\mathcal{A}' \otimes \mathcal{A}} (\widehat{CFDA}(\mathbb{I}), \widehat{CFDA}({\tau_{m+1}}) )$ with a direct system $(A_i, \phi_{ij} )_{i=0}^{\infty}$. Let $A_0$ be a type-$D$ module over $dg$ algebra $\mathcal{A}'$, and $A_i$ be a type-$D$ module homotopy equivalent to $A_0 \boxtimes (\underbrace{\tau_{m+1} \boxtimes \cdots \boxtimes \tau_{m+1}}_{i-\textrm{times}})$. Taking a tensor product with the identity type-$D$ morphism $\mathrm{id}$, we obtain $f \boxtimes \mathrm{id} : A_i \rightarrow A_{i+1}$. The map of the direct system is defined as
\begin{displaymath}
\phi_{i(i+j)} := \underbrace{ (f \boxtimes \mathrm{id}) \circ \cdots \circ (f \boxtimes \mathrm{id}) }_{j-\textrm{times} } : A_i \rightarrow A_{i+j}.
\end{displaymath}
Clearly, $\phi_{ii}$ is defined to be the identity morphism.
\begin{remark}
In general, the composition of type-$D$ morphism over $\mathcal{A}_{\infty}$-algebra is associative only up to homotopy. However, the torus algebra is $dg$ algebra, thus the composition is associative.
\end{remark}
The direct limit $A = \varinjlim A_i$, if exists, of the direct system $(A_i, \phi_{ij})$ has a differential induced from the differential of each type-$D$ module $A_i$. Equivalently, $\delta^1 [v] := [\delta^1 v]$. The algebra action is also defined to be $\rho_I \cdot [v] := [\rho_I \cdot v]$, where $I \in \{ 1,2,3,12,23,123 \}$ or $\rho_I$ is an idempotent (here, we view the type-$D$ structure as a $dg$  module). \\

In Section 5, we obtained four different nontrivial maps in $\textrm{Mor}_{\mathcal{A}' \otimes \mathcal{A}} (\widehat{CFDA}(\mathbb{I}), \widehat{CFDA}({\tau_{m+1}}) )$. However, there is only one element which gives an interesting direct system.
\begin{prop}
Let $f_i \in \textrm{Mor}_{\mathcal{A}' \otimes \mathcal{A}} (\widehat{CFDA}(\mathbb{I}), \widehat{CFDA}({\tau_{m+1}}) )$, $i=1, \cdots 4$. be maps in Section 5. For $f = a_1 \cdot f_1 + \cdots + a_4 \cdot f_4$, the direct system obtained by the above construction has trivial direct limit, unless $a_1 = 1$. 
\end{prop}
\emph{Proof.} Every value of the maps $f_2 \boxtimes \mathrm{id} , f_3 \boxtimes \mathrm{id}$ and $f_4 \boxtimes \mathrm{id}$ has a non-identity algebra element in $\mathcal{A}'$, and composing such maps will be eventually the zero map. $\square$ \\

From now on, we will be interested in the direct system constructed by $f_1 \boxtimes \mathrm{id}$, and $\phi_{ij}$ denotes the direct system map constructed from $f_1$.
\begin{prop}
Let $A_0$ be a type-$D$ module such that
\begin{displaymath}
\xymatrix{
\eta_{i+1} & \lambda^i_l \ar[l]_{\rho_2} & \cdots \ar[l]_{\rho_{23}} & \lambda^i_{k+1} \ar[l]_{\rho_{23}} & \lambda^i_k \ar[l]_{\rho_{23}} & \cdots \ar[l]_{\rho_{23}} & \lambda^i_1 \ar[l]_{\rho_{23}} & \eta_i \ar[l]_{\rho_3}.
}
\end{displaymath}
The direct limit of the direct system $(A_i, \phi_{ij})$ is isomorphic to $A_0$.
\label{prop:horizontal}
\end{prop}
To prove the above proposition, we will use maps described in Section 4. In order to do so, we introduce a little bit of an algebraic result.
\begin{lem}
Let $M$ be a type-$D$ module over $\mathcal{A}'$ of the following structure.
\begin{displaymath}
\xymatrix@R=1pc{
a \ar[d]^{\rho_1} & b \ar[l]_{\rho_2} & \ar@{..>}[l] \\
c & & \\
d \ar[u]_1 & &
}
\end{displaymath}
Then $M$ is homotopy equivalent to the following type-$D$ module.
\begin{displaymath}
\xymatrix{
x & y \ar[l]_{\rho_2} & \ar@{..>}[l] \\
}
\end{displaymath}
\label{lem:simplify3}
\end{lem}
\emph{Proof.} Name the latter type-$D$ module $M'$. We define two maps $F : M \rightarrow M'$ and $G : M' \rightarrow M$ such that
\begin{displaymath}
F(a) = x, \quad F(b) =y, \quad F(c)=F(d)=0,
\end{displaymath} 
and otherwise identity, and
\begin{displaymath}
G(x) = a + \rho_1 d, \quad G(y)=b
\end{displaymath}
and otherwise identity. Clearly $F \circ G =\mathbb{I}_{M'}$ and $G \circ F + \mathbb{I}_M = \delta^1 \circ H + H \circ \delta^1$ where $ H : M \rightarrow M$ is everywhere zero but $H(c)=d$. $\square$ \\

\begin{lem}
Let $M$ be a type-$D$ module that contains the following generators and differential.
\begin{displaymath}
\xymatrix{
\ar@{..>}[r]^{\rho_I} & w \ar[r]^{\rho_J} & z & y \ar[l]_1 \ar[r]^{\rho_K} & x \ar@{..>}[r]^{\rho_L} & ,
}
\end{displaymath}
where $\rho_J \cdot \rho_K = \rho_{JK} \neq 0$.
Then $M$ is homotopy equivalent to $M'$, where $M'$ is identical to $M$ but the above chain has been replaced to
\begin{displaymath}
\xymatrix{
\ar@{..>}[r]^{\rho_I} & b \ar[r]^{\rho_{JK}} & a \ar@{..>}[r]^{\rho_L} & .
}
\end{displaymath}
\label{lem:simplify1}
\end{lem}
\emph{Proof}. We define two type-$D$ module maps $F: M \rightarrow M'$ and $G : M' \rightarrow M$ as follows. First, $F$ is defined to be 
\begin{displaymath}
F(x) =a, \quad F(y) =0, \quad F(z)= \rho_K a, \quad \textrm{and} \quad F(w) =b
\end{displaymath}
and otherwise identity. In addition, $G$ is defined to be $G(b) = \rho_J y + w$ and $G(a) = x$, and otherwise identity, too. Then, $F \circ G + \mathbb{I}_{M'} = 0$ and $G \circ F + \mathbb{I}_M = \delta^1 \circ H + H \circ \delta^1$, where $H : M \rightarrow M$ is
\begin{displaymath}
H(x) = 0, \quad H(z) =y, \quad H(y)=0, \quad \textrm{and} \quad H(w)=0,
\end{displaymath}
and otherwise zero. $\square$\\
\begin{lem}
Suppose $\rho_J \cdot \rho_K = \rho_{JK} \neq 0$. Let $M_1$ be a type-$D$ module that contains the following generators and differentials.
\begin{displaymath}
\xymatrix{
& a \ar[d]^{\rho_{JK}} & b \ar[d]^{\rho_J} \ar[l]_1 \ar[r]^1 & c \ar[d]^{\rho_J} & d \ar[l]_1 \ar[d]^{\rho_J} & \\
& e \ar@{..>}[l]^{\rho_L} & f \ar[l]^{\rho_K} \ar[r]_1 & g & h \ar[l]^1 \ar@{..>}[r]_{\rho_I} & .
}
\end{displaymath}
Then $M_1$ is homotopy equivalent to $M_1'$, such that $M_1'$ has the same type-$D$ structure as $M$ but the above chain is replaced to
\begin{displaymath}
\xymatrix{
& x \ar@{..>}[l]^{\rho_L} & y \ar[l]^{\rho_K} \ar@{..>}[r]_{\rho_I} & .
}
\end{displaymath} 
In addition, let $M_2$ be a type-$D$ module that contains
\begin{displaymath}
\xymatrix{
\ar@{..>}[r]^{\rho_I} & a' \ar[d]^{\rho_J} & b' \ar[l]_1 \ar[d]^{\rho_J} \ar[r]^1 & c' \ar[d]^{\rho_{JK}} & \\
& d' & e' \ar[l]^1 \ar[r]_{\rho_K} & f' \ar@{..>}[r]_{\rho_L} &
}
\end{displaymath}
Then $M_2$ is homotopy equivalent to $M_2'$, where the above chain has been replaced to
\begin{displaymath}
\xymatrix{
\ar@{..>}[r]_{\rho_I}& z \ar[r]_{\rho_{JK}} & w \ar@{..>}[r]_{\rho_L} & .
}
\end{displaymath}
\label{lem:simplify2}
\end{lem}
\emph{Proof}. We show the chain complexes above can be projected to the chain complex $M$ in Lemma \ref{lem:simplify1}, without changing the homotopy type. For the first part of the proof, let $M_1''$ be the chain complex $M$ in Lemma \ref{lem:simplify1} where $\rho_J = 1$. Equivalently, $M_1''$ has the following structure.  
\begin{displaymath}
\xymatrix{
& w \ar@{..>}[l]_{\rho_I} \ar[r]^1 & z & y \ar[l]_1 \ar[r]^{\rho_K} & x \ar@{..>}[r]^{\rho_L} & .
}
\end{displaymath}
Then we use $F_1 : M_1 \rightarrow M_1''$ and $G_1 : M_1'' \rightarrow M_1$ such that
\begin{eqnarray*}
F_1 (a) = \rho_J (y + w), & F_1(e) = x, & F_1(f)=y, \\
F_1(g)=z, & F_1(c) = \rho_J w, & F_1(h) = w, \\
F_1(b) = 0, & F_1(d) =0, &
\end{eqnarray*}
and otherwise identity. Likewise,
\begin{displaymath}
G_1(x) =e, \quad G_1(y) =f, \quad G_1(z)=g, \quad G_1(w) =h,
\end{displaymath}
and otherwise identity. It is clear to show $F_1 \circ G_1 + \mathbb{I}_{M_1''} = 0$ and $G_1 \circ F_1 + \mathbb{I}_{M_1} = \delta^1 \circ H_1 + H_1 \circ \delta^1$, where $H_1 : M_1 \rightarrow M_1$ is
\begin{displaymath}
H(a) = b+ d, \quad H(c) = d
\end{displaymath}
and otherwise zero. Now it is a straightforward exercise to show $M_1'' \cong M_1'$. \\
The second part of the proof uses almost identical maps; let $M_2''$ be precisely the same complex as $M$ in Lemma \ref{lem:simplify1}. The maps we will use are $F_2 : M_2 \rightarrow M_2''$ and $G_2 : M_2'' \rightarrow M_2$ such that
\begin{eqnarray*}
F_2(a) = w, & F_2 (b) =0, & F_2(c) = \rho_J y + w, \\
F_2(d) = z, & F_2 (e) =y, & F_2(f) = x
\end{eqnarray*}
and 
\begin{displaymath}
G_2 (x) = f, \quad G_2(y) = e, \quad G_2(z) =d, \quad G_2(w) =a. \\
\end{displaymath}
(If not specified, everything else will be mapped to itself.) By defining a map $H_2 : M_2 \rightarrow M_2$ so that $H_2(c) =b$ and otherwise zero, it is again a straightforward exercise verifying $F_2 \circ G_2 + \mathbb{I}_{M_2''} =0$ and $G_2 \circ F_2 + \mathbb{I} = \delta^1 \circ H_2 + H_2 \circ \delta^1$, as well as showing $M_2'' \cong M_2'$. $\square$. \\

\emph{Proof of Proposition \ref{prop:horizontal}.} For computational convenience, we will use $ \widehat{CFDA}(\mathbb{I}') \boxtimes A_0$ as the model for $A_0$, and $\widehat{CFDA}(\tau_{m+1}) \boxtimes \widehat{CFDA}(\mathbb{I}') \boxtimes A_0 $ for $A_1$. The former complex has the following structure
\begin{displaymath}
\small
\xymatrix@C=1pc@R=1.5pc{
\mathbf{xx} \eta_{i+1} \ar@/_3pc/[dd]|{\rho_1} & \mathbf{yy} \lambda_l^i \ar[l]_{\rho_2} & \mathbf{yy} \lambda_{l-1}^i \ar[ddl]|{\rho_2} & & \mathbf{yy} \lambda_{k+1}^1 \ar[dl]|{\rho_2} & \mathbf{yy} \lambda_k^i \ar[ddl]|{\rho_2} & & \mathbf{yy} \lambda_1^i \ar[dl]|{\rho_2} & \mathbf{xx} \eta_i \ar@/^3pc/[dd]|{\rho_1} \ar[ddl] \\
\mathbf{yz_1} \eta_{i+1} \ar[d] & \mathbf{xw_1} \lambda_l^i \ar[u]_{\rho_3} \ar[d] & \mathbf{xw_1} \lambda_{l-1}^i \ar[u]_{\rho_3} \ar[d] & \cdots \ar[dl]|{\rho_2} & \mathbf{xw_1} \lambda_{k+1}^1 \ar[u]_{\rho_3} \ar[d] & \mathbf{xw_1} \lambda_k^i \ar[u]_{\rho_3} \ar[d] & \cdots \ar[dl]|{\rho_2} & \mathbf{xw_1} \lambda_1^i \ar[u]_{\rho_3} \ar[d] & \mathbf{yz_1} \eta_i \ar[d] \\
\mathbf{yz_2} \eta_{i+1} & \mathbf{xw_2} \lambda_l^i & \mathbf{xw_2} \lambda_{l-1}^i & & \mathbf{xw_2} \lambda_{k+1}^1 & \mathbf{xw_2} \lambda_k^i & & \mathbf{xw_2} \lambda_1^i & \mathbf{yz_2} \eta_i \\
}
\end{displaymath}
while the latter one has
\begin{displaymath}
\small
\xymatrix@C=1.25pc@R=1.5pc{
\mathbf{px} \eta_{i+1} \ar@/_2pc/[ddd]|{\rho_1} & \mathbf{qy} \lambda_l^i \ar[dl]|{\rho_{23}} & \mathbf{qy} \lambda_{l-1}^i \ar[dl]|{\rho_{23}} & & \mathbf{qy} \lambda_{k+1}^1 \ar[dl]|{\rho_{23}} & \mathbf{qy} \lambda_k^i \ar[dl]|{\rho_{23}} & & \mathbf{qy} \lambda_1^i \ar[dl]|{\rho_{23}} & \mathbf{px} \eta_i \ar@/^2pc/[ddd]|{\rho_1} \ar@/^/[ddl] \\
\mathbf{rx} \eta_{i+1} \ar[u]|{\rho_2} & \mathbf{rw_2} \lambda_l^i \ar[d]|{\rho_2} & \mathbf{rw_2} \lambda_{l-1}^i \ar[d]|{\rho_2} & & \mathbf{rw_2} \lambda_{k+1}^1 \ar[d]|{\rho_2} & \mathbf{rw_2} \lambda_k^i \ar[d]|{\rho_2} & & \mathbf{rw_2} \lambda_1^i \ar[d]|{\rho_2} & \mathbf{rx} \eta_i \ar[u]|{\rho_2} \ar[l] \\
\mathbf{qz_1} \eta_{i+1} \ar[d] & \mathbf{pw_2} \lambda_l^i & \mathbf{pw_2} \lambda_{l-1}^i & \cdots & \mathbf{pw_2} \lambda_{k+1}^1 & \mathbf{pw_2} \lambda_k^i & \cdots& \mathbf{pw_2} \lambda_1^i & \mathbf{qz_1} \eta_i \ar[d] \\
\mathbf{qz_2} \eta_{i+1} & \mathbf{rw_1} \lambda_l^i \ar@/^2pc/[uuu] \ar@/_1.5pc/[uu] \ar[d]|{\rho_2} & \mathbf{rw_1} \lambda_{l-1}^i \ar@/^2pc/[uuu] \ar@/_1.5pc/[uu] \ar[d]|{\rho_2} & & \mathbf{rw_1} \lambda_{k+1}^1 \ar@/^2pc/[uuu] \ar@/_1.5pc/[uu] \ar[d]|{\rho_2} & \mathbf{rw_1} \lambda_k^i \ar@/^2pc/[uuu] \ar@/_1.5pc/[uu] \ar[d]|{\rho_2} & & \mathbf{rw_1} \lambda_1^i \ar@/^2pc/[uuu] \ar@/_1.5pc/[uu] \ar[d]|{\rho_2} & \mathbf{qz_2} \eta_i \\
& \mathbf{pw_1} \lambda_l^i \ar@/^1.5pc/[uu] \ar@/^0.75pc/[uuul]|{\rho_3} & \mathbf{pw_1} \lambda_{l-1}^i \ar@/^1.5pc/[uu] \ar@/^0.75pc/[uuul]|{\rho_3} & & \mathbf{pw_1} \lambda_{k+1}^1 \ar@/^1.5pc/[uu] \ar@/^0.75pc/[uuul]|{\rho_3} & \mathbf{pw_1} \lambda_k^i \ar@/^1.5pc/[uu] \ar@/^0.75pc/[uuul]|{\rho_3} & & \mathbf{pw_1} \lambda_1^i. \ar@/^1.5pc/[uu] \ar@/^0.75pc/[uuul]|{\rho_3} & \\
}
\end{displaymath}
Note that these two complexes are homotopy equivalent to $A_0$.\\

By Lemma~\ref{lem:simplify3}, $\mathbf{yz_1} \eta_{i+1} \rightarrow \mathbf{yz_2} \eta_{i+1}$ and $\mathbf{yz_1} \eta_i \rightarrow \mathbf{yz_2} \eta_i$ from $ \widehat{CFDA}(\mathbb{I}') \boxtimes A_0$, and $\mathbf{qz_1} \eta_{i+1} \rightarrow \mathbf{qz_2} \eta_{i+1}$ and $\mathbf{qz_1} \eta_i \rightarrow \mathbf{qz_2} \eta_i$ from $\widehat{CFDA}(\tau_{m+1}) \boxtimes \widehat{CFDA}(\mathbb{I}') \boxtimes A_0 $ can be removed without changing the homotopy type. \\

We now take an identity map $\mathrm{id}:=\mathrm{id}_{A_0}$ and take a tensor product with the map $(f_1 \boxtimes \mathrm{id})$. The computation of the resulting map $f_1 \boxtimes \mathrm{id} \cong \phi_{01}$  is straightforward by the Pairing Theorem, which is from $ \widehat{CFDA}(\mathbb{I}') \boxtimes A_0$ to $\widehat{CFDA}(\tau_{m+1}) \boxtimes \widehat{CFDA}(\mathbb{I}') \boxtimes A_0 $. However we can greatly simplify the map to the map from $A_0$ to $A_1 \cong A_0$, by repeatedly applying Lemma \ref{lem:simplify1} and Lemma \ref{lem:simplify2}. 
\begin{itemize}
  \item $\mathbf{yy} \lambda_l^i \ \mapsto \ \mathbf{rx} \eta_{i+1}$. The generator $\mathbf{yy} \lambda_l^i$ corresponds to $x \in M$ in Lemma \ref{lem:simplify1}, and $\mathbf{rx} \eta_{i+1}$ corresponds to $f' \in M_2$ in Lemma \ref{lem:simplify3}. This implies the image of $\lambda_l^i$ under $\phi_{01}$ has a term of $\lambda_l^i$.  
  \item $\mathbf{yy} \lambda_m^i \ \mapsto \ \mathbf{rw_2} \lambda_{m+1}^i$, where $m = 1, \cdots, l-1$. Similarly as above, this can be interpreted as the image of $\lambda_m^i$ having a term of $\lambda_m^i$. 
  \item $\mathbf{xw_j} \lambda_m^i \ \mapsto \ \mathbf{pw_j} \lambda_m^i + \rho_3 \cdot \mathbf{rw_j} \lambda_m^i$, where $m = 1, \cdots l$ and $j=1, 2$.	Due to Lemma \ref{lem:simplify2}, $\mathbf{pw_j} \lambda_m^i$ and  $\mathbf{rw_j} \lambda_m^i$ are homotopic to zero, thus it has no contribution towards the map $\phi_{01}$. 
  \item $\mathbf{xx} \eta_j \ \mapsto \ \mathbf{px} \eta_j + \rho_3 \cdot \mathbf{rx} \eta_j$, where $j=i, i+1$. First we consider the case when $j=i$. Note $\mathbf{xx} \eta_i \cong \eta_i \in A_0$. Again by Lemma \ref{lem:simplify2}, $\mathbf{rx} \eta_i$ is homotopic to zero and $\mathbf{px} \eta_i \cong \eta_i \in A_1$, thus it can be interpreted as $ \eta_i \mapsto \eta_i$. Second, if $j=i+1$, $\mathbf{px} \eta_{i+1} \cong \eta_{i+1} \in A_1$ and $\mathbf{rx} \eta_{i+1} \cong \lambda_l^i \in A_1$, which means $\eta_{i+1} \mapsto \eta_{i+1} + \rho_3 \cdot \lambda_l^i$.
  \item $\mathbf{yy} \lambda_m^i \ \mapsto \ \mathbf{qy} \lambda_m^i$, where $m=1, \cdots, l$. It is easy to show $\mathbf{yy} \lambda_m^i \cong \lambda_m^i \in A_0$, and $\mathbf{qy} \lambda_m^i \cong \lambda_{m-1} \in A_1$ for $m=2, \cdots l$, by again applying Lemma \ref{lem:simplify2}. For $m=1$ case, $\mathbf{yy} \lambda_1^i$ is still homotopic to $\lambda_l^1 \in A_0$ but $\mathbf{qy} \lambda_1^i \cong \rho_2 \cdot \eta_i \in A_1$ by the first part of Lemma \ref{lem:simplify2}. 
\end{itemize}
Summarizing the above result, the map $\phi_{01}$ is
\begin{eqnarray*}
\phi_{01} (\eta_i) & = & \eta_i \\
\phi_{01} (\eta_{i+1}) & = & \eta_{i+1} + \rho_3 \cdot \lambda_l^i \\
\phi_{01} (\lambda_m^i) & = & \lambda_m^i + \lambda_{m-1}^i, \textrm{ where } m = 2, \cdots l. \\
\phi_{01} (\lambda_1^i) & = & \lambda_1^i + \rho_2 \cdot \eta_i.
\end{eqnarray*}
Note that the proof so far only concerned the map $\phi_{01} : A_0 \rightarrow A_1$, but the behavior of maps $\phi_{i,i+1}$ is exactly the same as $\phi_{01}$. See Figure \ref{figure:phi01}. \\ 

\begin{figure}
\begin{displaymath}
\xymatrix@R=3pc{
A_i : \ar[d]|{\phi_{i,i+1}} & \eta_{i+1} \ar[d] \ar[dr]|{\rho_3} & \lambda^i_l \ar[l]_{\rho_2} \ar[d] \ar[dr] & \cdots \ar[l]_{\rho_{23}} & \lambda^i_{k+1} \ar[l]_{\rho_{23}} \ar[d] \ar[dr] & \lambda^i_k \ar[l]_{\rho_{23}} \ar[d] \ar[dr] & \cdots \ar[l]_{\rho_{23}} & \lambda^i_1 \ar[l]_{\rho_{23}} \ar[d] \ar[dr]|{\rho_2} & \eta_i \ar[l]_{\rho_3} \ar[d] \\
A_{i+1} : & \eta_{i+1} & \lambda^i_l \ar[l]_{\rho_2} & \cdots \ar[l]_{\rho_{23}} & \lambda^i_{k+1} \ar[l]_{\rho_{23}} & \lambda^i_k \ar[l]_{\rho_{23}} & \cdots \ar[l]_{\rho_{23}} & \lambda^i_1 \ar[l]_{\rho_{23}} & \eta_i \ar[l]_{\rho_3}
}
\end{displaymath}
\caption{The above row represents $A_i$ and the below $A_{i+1}$, and the arrows from the top row to the bottom row represents the map $\phi_{i,i+1}$. }
\label{figure:phi01}
\end{figure}
We now turn to the direct system. First, note the map can be represented by a matrix similar to the Jordan matrix. In particular, the ordered basis $\{ \eta_1, \lambda_1^i, \cdots, \lambda_m^i, \eta_{i+1} \}$ gives following matrix.
\begin{displaymath}
X:=[\phi_{i,i+1}] =
\left( \begin{array}{cccccc}
  	1 & \rho_2 &   &        &   & \\
  	  & 1      & 1 &        &   & \\
  	  &        & 1 &        &   & \\
  	  &        &   & \ddots &   & \\
  	  &        &   &        & 1 & \rho_3 \\
  	  &        &   &        &   & 1 
\end{array}
\right)
\end{displaymath}  
First, note that $\phi_{i,j} \cong \phi_{i,i+1} \circ \phi_{i+1,i+2} \circ \cdots \phi_{j-1,j} $. Since the coefficients appearing in this matrix is in $\mathbb{F}_2$, there is an integer $\nu=\nu(l)$, which depends on the length of the chain $A_i$ such that $X^{\nu}$ is the identity matrix. Thus, we can choose a subsystem $A_i, A_{i+\nu}, A_{i+2 \nu}, \cdots$ with every map in the subsystem is identity. By using cofinality the proposition is proved. $\square$ \\

By parallel computation, the following proposition regarding vertical complex is easily proved.
\begin{prop}
Let $A_0$ be a type-$D$ module such that 
\begin{displaymath}
\xymatrix{
\xi_i \ar[r]^{\rho_1} & \kappa^i_1 & \cdots \ar[l]_{\rho_{23}} & \kappa^i_k \ar[l]_{\rho_{23}} & \kappa^i_{k+1} \ar[l]_{\rho_{23}} & \cdots \ar[l]_{\rho_{23}} & \kappa^i_l \ar[l]_{\rho_{23}} & \xi_{i+1} \ar[l]_{\rho_{123}}.
}
\end{displaymath}
Then the direct limit of the direct system $(A_i, \phi_{ij})$ is isomorphic to $A_0$ as well. 
\label{prop:vertical}
\end{prop}
\emph{Proof.} The proof is similar; the map $\phi_{i,i+1} : A_i \rightarrow A_{i+1}$ is obtained by exactly the same manner as above, and it is given in the diagram below.
\begin{displaymath}
\xymatrix@C=2pc@R=3pc{
A_i :& \xi_i \ar[r]^{\rho_1} \ar[d] & \kappa^i_1 \ar[d] & \cdots \ar[l]_{\rho_{23}} \ar[dl] & \kappa^i_k \ar[l]_{\rho_{23}} \ar[d] \ar[dl] & \kappa^i_{k+1} \ar[l]_{\rho_{23}} \ar[d] \ar[dl] & \cdots \ar[l]_{\rho_{23}} \ar[dl] & \kappa^i_l \ar[l]_{\rho_{23}} \ar[dl] \ar[d] & \xi_{i+1} \ar[l]_{\rho_{123}} \ar[d] \ar[dl]|{\rho_1} \\
A_{i+1} :& \xi_i \ar[r]^{\rho_1} & \kappa^i_1 & \cdots \ar[l]_{\rho_{23}} & \kappa^i_k \ar[l]_{\rho_{23}} & \kappa^i_{k+1} \ar[l]_{\rho_{23}} & \cdots \ar[l]_{\rho_{23}} & \kappa^i_l \ar[l]_{\rho_{23}} & \xi_{i+1}. \ar[l]_{\rho_{123}}
}
\end{displaymath}
Then the proposition is justified by using the cofinality of the direct system. $\square$\\

\begin{remark}
Although Propositions \ref{prop:horizontal} and \ref{prop:vertical} show the structure of the direct limit of two different chain complexes, this proof can also be applied to the combination of these two chains. For example, the map $\phi_{ij}$ maps the following generators
\begin{displaymath}
\xymatrix@R=1pc@C=1pc{
\eta_1 \ar[d]^{\rho_1} & \lambda^1_l \ar[l]^{\rho_2} & \cdots \ar[l]^{\rho_{23}} & \lambda^1_{k+1} \ar[l]^{\rho_{23}} & \lambda^1_k \ar[l]^{\rho_{23}} & \cdots \ar[l]^{\rho_{23}} & \lambda^1_1 \ar[l]^{\rho_{23}} & \eta_2 \ar[l]^{\rho_3} \ar[d]^{\rho_1} \\
\kappa^1_1 & & & & & & & \kappa^2_1 \\
\vdots \ar[u]_{\rho_{23}} & & & & & & & \vdots \ar[u]_{\rho_{23}}\\
\kappa^1_k \ar[u]_{\rho_{23}}& & & & & & & \kappa^2_k \ar[u]_{\rho_{23}}\\
\kappa^1_{k+1} \ar[u]_{\rho_{23}}& & & & & & & \kappa^2_{k+1} \ar[u]_{\rho_{23}} \\
\vdots \ar[u]_{\rho_{23}} & & & & & & & \vdots \ar[u]_{\rho_{23}} \\
\kappa^1_m \ar[u]_{\rho_{23}} & & & & & & & \kappa^2_m \ar[u]_{\rho_{23}} \\
\eta_3 \ar[u]_{\rho_{123}} & \lambda^2_l \ar[l]^{\rho_2} & \cdots \ar[l]^{\rho_{23}} & \lambda^2_{k+1} \ar[l]^{\rho_{23}} & \lambda^2_k \ar[l]^{\rho_{23}} & \cdots \ar[l]^{\rho_{23}} & \lambda^2_1 \ar[l]^{\rho_{23}} & \eta_4 \ar[l]^{\rho_3} \ar[u]_{\rho_{123}}
}
\end{displaymath}
such that $\lambda^i_{k+1} \mapsto \lambda^i_{k+1} + \lambda^i_k$, $\kappa^i_{k+1} \mapsto \kappa^i_{k+1} + \kappa^i_k$ and so on. Especially we emphasize the images of $\eta_j$s. They are mapped to
\begin{itemize}
  \item $\eta_1 \mapsto \eta_1 + \rho_3 \cdot \lambda^1_l$
  \item $\eta_2 \mapsto \eta_2$
  \item $\eta_3 \mapsto \eta_3 + \rho_1 \cdot \kappa^1_m + \rho_3 \cdot \lambda^2_l$
  \item $\eta_4 \mapsto \eta_4 + \rho_1 \cdot \kappa^2_m$.
\end{itemize}
\end{remark}

Last but not least, we study the direct limit of unstable chain.

\begin{figure}
\begin{displaymath}
\xymatrix{
\xi_0^0 \ar[d] \ar[dr]|{\rho_1} & \mu_1^0 \ar[d] \ar[dr] & \cdots & \mu_l^0 \ar[d] \ar[dr] & \eta_0^0 \ar[d]|{\rho_3} \ar[dr] & & & & \\
\xi_0^1 \ar[d] \ar[dr]|{\rho_1} & \mu_1^1 \ar[d] \ar[dr] & \cdots & \mu_l^1 \ar[d] \ar[dr] & \mu_{l+1}^1 \ar[d] \ar[dr] & \eta_0^1 \ar[d]|{\rho_3} \ar[dr]& & & \\
& & & \vdots & & & \ddots & & \\
*+[F]{\xi_0^m} \ar[d] \ar[dr]|{\rho_1} & *+[F]{\mu_1^m} \ar[d] \ar[dr] & \cdots & *+[F]{\mu_{l+m-2}^m} \ar[d] \ar[dr] & *+[F]{\mu_{l+m-1}^m} \ar[d] \ar[dr] & *+[F]{\mu_{l+m}^m} \ar[d] \ar[dr] & *+[F]{\eta_0^m} \ar[d]|{\rho_3} \ar[dr] & & \\
\xi_0^{m+1} \ar[d] \ar[dr]|{\rho_1} & \mu_1^{m+1} \ar[d] \ar[dr] & \cdots & \mu_{l+m-2}^{m+1} \ar[d] \ar[dr] & \mu_{l+m-1}^{m+1} \ar[d] \ar[dr] & \mu_{l+m}^{m+1} \ar[d] \ar[dr] & *+[F]{\mu_{l+m+1}^{m+1}} \ar[d] \ar[dr] & \eta_0^{m+1} \ar[d]|{\rho_3} \ar[dr] & \\
& & \vdots & & & & & *+[F]{ \mu^{m+2}_{l+m+2} } & \ddots
}
\end{displaymath}
\caption{Each arrow represents a term of function $\phi_{i,i+1}$. The generators in the boxes are chosen as basis of the resulting direct system.}
\label{fig:dsystem}
\end{figure}
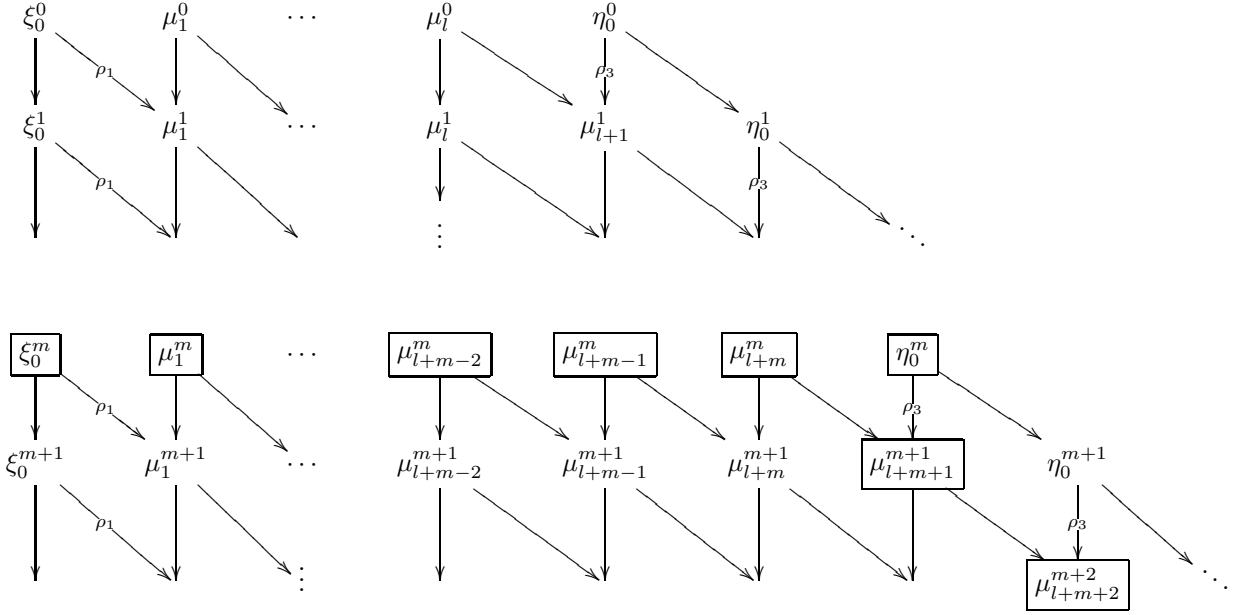

\begin{prop}
Suppose $A_0$ has the following structure. 
\begin{displaymath}
\xymatrix{
\xi_0^0 \ar[r]^{\rho_{123}} & \mu_1^0 \ar[r]^{\rho_{23}} & \mu_2^0 \ar[r]^{\rho_{23}} & \cdots \ar[r]^{\rho_{23}} & \mu_l^0 \ar[r]^{\rho_2} & \eta_0^0.
}
\end{displaymath} 
Then the direct limit $A = \varinjlim A_i$ of the direct system $(A_i, \phi_{ij})$ has the following structure. Let $m$ be an arbitrary integer greater than $l$. There are generators $\overline{\xi}$, $\overline{\eta}$ and $\overline{\mu}_1, \cdots, \overline{\mu}_m$ equipped with type-$D$ structure as below.
\begin{displaymath}
\xymatrix{
\overline{\xi} \ar[r]^{\rho_{123}} & \overline{\mu}_1 \ar[r]^{\rho_{23}} & \overline{\mu}_2 \ar[r]^{\rho_{23}} & \cdots \ar[r]^{\rho_{23}} & \overline{\mu}_m \ar[r]^{\rho_2} & \overline{\eta}.
}
\end{displaymath}
In addition, there are a set of generators $\{ \nu_i \}_{i=1}^{\infty}$ such that
\begin{displaymath}
\delta^1(\nu_i) = \sum_{k=1}^i \rho_{23} \cdot \nu_k + \rho_2 \cdot \overline{\eta}.
\end{displaymath} 
\end{prop}
\emph{Proof.} We construct the direct system $(A_i, \phi_{ij})$ as in Proposition \ref{prop:horizontal}, and the relation between $A_i$ and $A_{i+1}$ is given below.
\begin{displaymath}
\xymatrix{
A_i : \ar[d]|{\phi_{i,i+1}} & \xi_0^i \ar[r]^{\rho_{123}} \ar[d] \ar[dr]|{\rho_1} & \mu_1^i \ar[r]^{\rho_{23}} \ar[d] \ar[dr] & \cdots \ar[r]^{\rho_{23}} \ar[dr] & \mu_j^i \ar[r]^{\rho_{23}} \ar[d] \ar[dr] & \cdots \ar[r]^{\rho_{23}} & \mu_i^i \ar[r]^{\rho_2} \ar[d] \ar[dr] & \eta_0^i \ar[d]|{\rho_3} \ar[dr] & \\
A_{i+1} : & \xi_0^{i+1} \ar[r]^{\rho_{123}} & \mu_1^{i+1} \ar[r]^{\rho_{23}} & \cdots \ar[r]^{\rho_{23}} & \mu_j^{i+1}  \ar[r]^{\rho_{23}} & \mu_{j+1}^{i+1} \ar[r]^{\rho_{23}} & \cdots \ar[r]^{\rho_{23}} & \mu_{i+1}^{i+1} \ar[r]^{\rho_2} & \eta_0^{i+1}.
}
\end{displaymath}
See figure~\ref{fig:dsystem} for the diagram of the entire direct system. \\

The resulting direct system $\varinjlim A_i := \bigcup_{i=0}^{\infty} A_i / \sim$, where $x \sim y$ iff $x=\phi_{ij}(y)$. As the generators of the $\varinjlim A_i$, we choose $\xi_0^m, \eta_0^m, \mu_1^m, \cdots \mu_{l+m}^m \in A_m$ and $\mu_{l+m+k}^{m+k} \in A_{m+k}$. for $k =1 ,2, \cdots$. Now it is clear that every element in $\varinjlim A_i$ can be written uniquely as sum of chosen elements, possibly with the torus algebra action. We let $\overline{\xi}:= \xi_0^m$, $\overline{\eta} := \eta_0^m$, $\overline{\mu}_i := \mu_i^m$, and $\nu_i := \mu_{l+m+i}^{m+i}$. Then $\overline{\xi}$, $\overline{\mu}_i$ and $\overline{\eta}$ will clear have the type-$D$ structure as stated in the proposition. The differential of $\nu_i$ is also easily proved, since $\delta^1 (\mu_{l+m+i}^{m+i}) = \rho_2 \cdot \eta_0^{m+i} $ and $\eta_0^{m+i} = \phi_{m+i-1, m+i} (\eta_0^{m+i-1}) + \rho_3 \cdot \mu_{l+m+i}^{m+i}$. $\square$ \\

\textbf{Examples} The above propositions discussed about the separate chains which can be graded by the grading set $G'(T^2)$. Here we compute examples that cannot be graded by $G'(T^2)$, the right handed trefoil and the figure eight knot. The reduced type-$D$ module of the knot complement of the left handed trefoil has two distinguished generators $\eta_0$ and $\xi_0$, where $\eta_0 \neq \xi_0$. On the other hand, the figure eight knot has $\eta_0 = \xi_0$. \\

First we compute the figure eight knot. Let $A_0$ be the reduced type-$D$ module $\widehat{CFD}(S^3 \backslash \textrm{nbd}(K))$ of the figure eight knot $K$ complement with sufficiently large positive framing. The type-$D$ structure can be derived from the knot Floer homology $CFK^-(K)$, or from the Figure 11.5 of \cite{LOT08}. Then the type-$D$ structure of $\varinjlim A_i$ is depicted as below.
\begin{displaymath}
\xymatrix{
\eta_1 \ar[d]^{\rho_1} & \lambda^1_1 \ar[l]^{\rho_2} & \eta_2 \ar[l]^{\rho_3}  \ar[d]^{\rho_1} & & & \cdots \ar@/^1pc/[dr]|{\rho_{23}} & \\
\kappa^1_1 & & \kappa^2_1 & & \overline{\mu}_1 \ar@/^1pc/[ur]|{\rho_{23}} & & \overline{\mu}_m \ar@/^1pc/[dl]|{\rho_2} \\
\eta_3 \ar[u]_{\rho_{123}} & \lambda^2_1 \ar[l]^{\rho_2} & \eta_4 \ar[l]^{\rho_3} \ar[u]_{\rho_{123}} & & & \overline{\xi} \ar@/^1pc/[ul]|{\rho_{123}} &  \\
& & & & & \nu_1 \ar@(ul,dl)[]|{\rho_{23}} \ar[u]|{\rho_2} & \\
& & & & & \nu_2 \ar@(ul,dl)[]|{\rho_{23}} \ar[u]|{\rho_{23}} \ar@/_1.5pc/[uu]|{\rho_2} & \\
& & & & & \vdots &
}
\end{displaymath}
\bigskip

Second, we consider the case where $K=$right handed trefoil knot. Let $A_0$ be $\widehat{CFD}(S^3 \backslash \textrm{nbd}(K))$ with sufficiently large positive framing, again. Then the direct limit $\varinjlim A_i$ has following type-$D$ module.
\begin{displaymath}
\xymatrix@C=0.75pc@R=0.75pc{
\overline{\xi} \ar[dr]_{\rho_{123}} & & \lambda_1 \ar[ll]_{\rho_2} & & \eta_2 \ar[ll]_{\rho_3} \ar[dd]^{\rho_1} \\
& \overline{\mu}_1 \ar[dr]_{\rho_{23}} & & & \\
& & \ddots \ar[dr]_{\rho_{23}} & & \lambda_2 \\
& & & \overline{\mu}_m \ar[dr]_{\rho_2} & \\
\cdots \ar[r] & \nu_2 \ar@(dl,dr)[]|{\rho_{23}} \ar[rr]_{\rho_{23}} \ar@/_2.5pc/[rrr]_{\rho_2}& & \nu_1 \ar@(dl,dr)[]|{\rho_{23}} \ar[r]_{\rho_2} & \overline{\eta} \ar[uu]_{\rho_{123}}
}
\end{displaymath} 

\section{Proof of Theorem 2}
So far, we considered the knot complement equipped with sufficiently large positive framings. However even when the knot complement is equipped with sufficiently large negative framings, the behaviour of the map $\phi_{i,i+1}$ is easily understood. \\

For a type-$D$ module of knot $K \subset S^3$ complement with sufficiently large negative framing $n$, there are chains between the two distinguished generators $\eta_0$ and $\xi_0$ as follows.
\begin{displaymath}
\xymatrix{
\xi_0 \ar[r]|{\rho_1} & \mu_1 & \mu_2 \ar[l]|{\rho_{23}} & \cdots \ar[l]|{\rho_{23}} & \mu_m \ar[l]|{\rho_{23}} & \eta_0, \ar[l]|{\rho_3} 
}
\end{displaymath}
where $m=2 \tau(K) + n$.  

\begin{prop}
Let $A_i$ be a type-$D$ module such that 
\begin{displaymath}
\xymatrix{
\xi^i_0 \ar[r]|{\rho_1} & \mu^i_1 & \mu^i_2 \ar[l]|{\rho_{23}} & \cdots \ar[l]|{\rho_{23}} & \mu^i_{l-1} \ar[l]|{\rho_{23}} & \mu^i_l \ar[l]|{\rho_{23}} & \eta^i_0. \ar[l]|{\rho_3} 
}
\end{displaymath}
Then the map $\phi_{i,i+1} : A_i \rightarrow A_{i+1}$ is as follows.
\begin{displaymath}
\xymatrix{
A_i : & \xi^i_0 \ar[d] \ar[r]|{\rho_1} & \mu^i_1 \ar[d] & \mu^i_2 \ar[l]|{\rho_{23}} \ar[d] \ar[dl] & \cdots \ar[l]|{\rho_{23}} \ar[d] \ar[dl] & \mu^i_{l-1} \ar[l]|{\rho_{23}} \ar[d] \ar[dl] & \mu^i_l \ar[l]|{\rho_{23}} \ar[d]|{\rho_2} \ar[dl] & \eta^i_0 \ar[dl] \ar[l]|{\rho_3} \\
A_{i+1} : & \xi^{i+1}_0 \ar[r]|{\rho_1} & \mu^{i+1}_1 & \mu^{i+1}_2 \ar[l]|{\rho_{23}} & \cdots \ar[l]|{\rho_{23}} & \mu^{i+1}_{l-1} \ar[l]|{\rho_{23}} & \eta^{i+1}_0 \ar[l]|{\rho_3} &
}
\end{displaymath}
\end{prop}
\emph{Proof.} Straightforward. $\square$ \\

\begin{figure}
\begin{displaymath}
\xymatrix@R=1.7pc{
\xi^0_0 \ar[d] \ar[r]|{\rho_1} & \mu^0_1 \ar[d] & \mu^0_2 \ar[l]|{\rho_{23}} \ar[d] \ar[dl] & \cdots \ar[l]|{\rho_{23}} \ar[d] \ar[dl] & \mu^0_{l-1} \ar[l]|{\rho_{23}} \ar[d] \ar[dl] & \mu^0_l \ar[l]|{\rho_{23}} \ar[d]|{\rho_2} \ar[dl] & \eta^0_0 \ar[dl] \ar[l]|{\rho_3} \\
\xi^1_0 \ar[d] \ar[r]|{\rho_1} & \mu^1_1 \ar[d] & \mu^1_2 \ar[l]|{\rho_{23}} \ar[d] \ar[dl] & \cdots \ar[l]|{\rho_{23}} \ar[d] \ar[dl] & \mu^1_{l-1} \ar[l]|{\rho_{23}} \ar[d]|{\rho_2} \ar[dl] & \eta^1_0 \ar[dl] \ar[l]|{\rho_3} & \\
& & & & & & \\
& \vdots & & \vdots & & & \\
\xi^{l-2}_0 \ar[d] \ar[r]|{\rho_1} & \mu^{l-2}_1 \ar[d] & \mu^{l-2}_2 \ar[l]|{\rho_{23}} \ar[d]|{\rho_2} \ar[dl] & \eta^{l-2}_0 \ar[dl] \ar[l]|{\rho_3} & & & \\
\xi^{l-1}_0 \ar[r]|{\rho_1} \ar[d] & \mu^{l-1}_1 \ar[d]|{\rho_2} & \eta^{l-1}_0 \ar[l]|{\rho_3} \ar[dl] & & & & \\
\xi^l_0 \ar[r]|{\rho_{12}} \ar[d] \ar[dr]|{\rho_1} & \eta^l_0 \ar[d]|{\rho_3} \ar[dr] & & & & & \\
\xi^{l+1}_0 \ar[r]|{\rho_{123}} \ar[d] \ar[dr]|{\rho_1} & \mu^{l+1}_1 \ar[r]|{\rho_2} \ar[d] \ar[dr] & \eta^{l+1}_0 \ar[d]|{\rho_3} \ar[dr] & & & & \\
\xi^{l+2}_0 \ar[r]|{\rho_{123}} & \mu^{l+2}_1 \ar[r]|{\rho_{23}} & \mu^{l+2}_2 \ar[r]|{\rho_2} & \eta^{l+2}_0 & & & \\
& & & \vdots & & &
}
\end{displaymath}
\caption{The direct system $(A_i, \phi_{ij})_{i=0}$, with $A_0$ equal to the unstable chain of sufficiently large negative framing. The image of $A_0$ under $\phi_{0,l}$ is homotopic to the complex that consists of two generators such that $\xi \mapsto \rho_{12} \cdot \eta$. In particular, for the natural inclusion map $i : A_0 \rightarrow A_l$, $i(\iota_0 \cdot A_0 ) \cong 0$. }
\label{fig:unstable}
\end{figure}
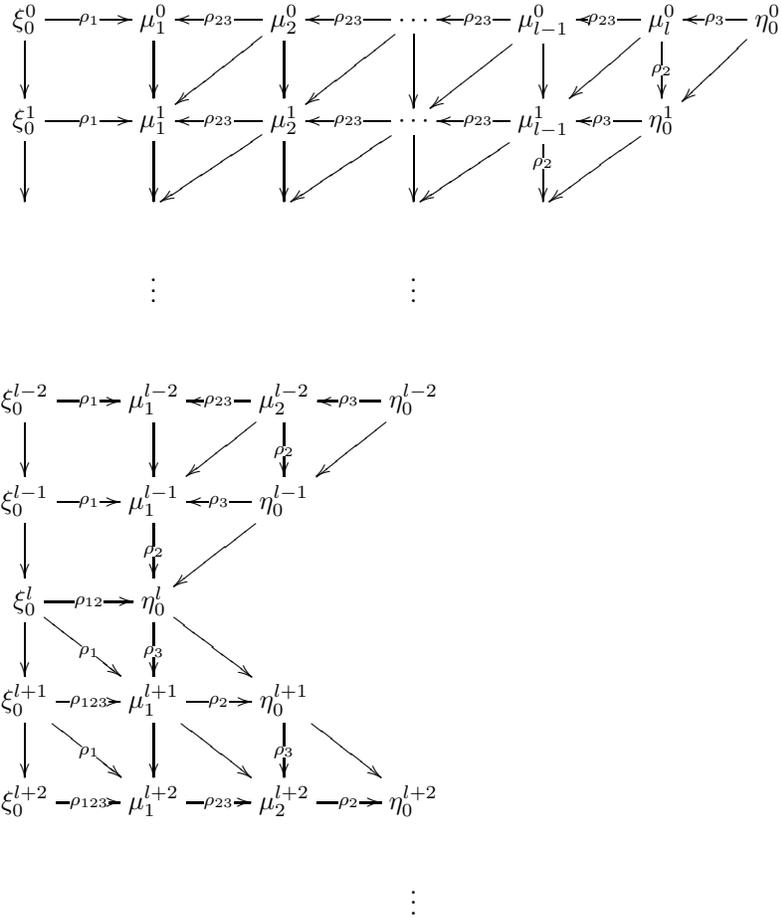

\emph{Proof of the Theorem 2.} We now consider a direct system $(A_i, \phi_{ij})_{i=0}$, where $A_0$ is a type-$D$ module of a knot $K$ complement with sufficiently large negative framing $i_0$. Also assume $A_i$ is the chain complex obtained by using the algorithm in Theorem 11.26 in~\cite{LOT08}. For negative $i$ values, the unstable chain of $A_i$ eventually become equal, as the framing number passes from negative to positive integer. In other words, 
\begin{displaymath}
\phi_{0, |2 \tau(K) - i_0|} : A_0 \rightarrow A_{|2 \tau(K) - i_0|}
\end{displaymath}
maps all unstable chains in $A_0$ to zero (see the figure~\ref{fig:unstable}). Thus the natural inclusion $\iota(A_0)$ is a knot invariant. $\square$

\end{document}